\documentclass{amsart}
\usepackage{amssymb}
\usepackage{mathabx}	

\usepackage[dvips]{graphicx,color}
\usepackage[percent]{overpic}
\usepackage[all]{xy}
\usepackage{tabularx}

\newcommand{\CC}{{\mathcal C}}          
\newcommand{\DD}{{\mathcal D}}		
\newcommand{\V}{\mathcal V} 		
\newcommand{\Ve}{\widetilde{\mathcal V}} 
\newcommand{\Rho}{R}			
\renewcommand{\1}{\mathbf{1}}		
\newcommand{\e}{1}   		 	
\newcommand{\U}{\mathcal{U}}		
\newcommand{\basis}[2][g]{\vphantom{\lambda}^{#1}{\lambda}_{#2}}
\newcommand{\ph}{\varphi}
\newcommand{\C}{\mathbb{C}}       	
\newcommand{\de}{\delta}
\newcommand{\xxto}{\xrightarrow}        

\newcommand{\ttt}{\otimes}              
\newcommand{\sz}{\scriptsize}		
\newcommand{\fz}{\footnotesize}         
\newcommand{\Z}{\mathbb{Z}}		
\newcommand{\isoto}{\widetilde{\rightarrow}}
\newcommand{\isomor}{\widetilde{-\mspace{-11mu}-\mspace{-11mu}-}}
\newcommand{\injto}{\hookrightarrow}    

\newcommand{\sem}{\sf}			
\newcommand{\dem}{\sf}			

\newcommand{\paref}[1]{Part~{\rm\ref{#1}}}
\newcommand{\firef}[1]{Figure~{\rm\ref{#1}}}
\newcommand{\thref}[1]{Theorem~{\rm\ref{#1}}}

\newcommand{\leref}[1]{Lemma~{\rm\ref{#1}}}
\newcommand{\coref}[1]{Corollary~{\rm\ref{#1}}}
\newcommand{\deref}[1]{Definition~{\rm\ref{#1}}}

\newcommand{\seref}[1]{Section~{\rm\ref{#1}}}

\swapnumbers 		
\newtheorem{theorem}[subsection]{Theorem}
\newtheorem{lemma}[subsection]{Lemma}
\newtheorem{proposition}[subsection]{Proposition}
\newtheorem{corollary}[subsection]{Corollary}

\theoremstyle{definition}
\newtheorem{definition}[subsection]{Definition}
\newtheorem{notation}[subsection]{Notation}

\theoremstyle{remark}
\newtheorem{nb}[subsection]{Nota Bene}

\DeclareMathOperator{\Mor}{Mor}
\DeclareMathOperator{\id}{id}

\DeclareMathOperator{\rep}{Rep}
\DeclareMathOperator{\Hom}{Hom}

\begin{document}

\title[Fusion in the Extended Verlinde Algebra]{Fusion in the Extended Verlinde Algebra}
\author{Vincent Graziano JR}

\address{DEPT of MATHS, SUNY at Stony Brook, Stony Brook, NY 11794, USA}
\email{graziano@math.sunysb.edu}
\urladdr{http://www.math.sunysb.edu/\textasciitilde graziano/}
\date{\today}

\maketitle
\thispagestyle{empty} 

\setcounter{tocdepth}{1}

\section*{Introduction} 

In this paper we present results which allow one to work with $G$-equivariant fusion categories. These categories were formalized so as to capture many of the features of conformal field theories that have a finite group of automorphisms acting on them. Here, in some sense, we bear the first fruits of this formalization.

\begin{figure}[ht]
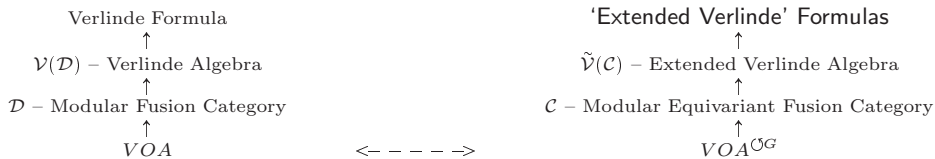

 \begin{scriptsize}
  \begin{tabularx}{\linewidth}{cXcXc}
    Verlinde Formula &&&&{\sem \small{`Extended Verlinde' Formulas}}\\
  $\uparrow$ &&&& $\uparrow$ \\
$\V(\DD)$  -- Verlinde Algebra &&&& $\Ve(\CC)$ -- Extended Verlinde Algebra\\
$\uparrow$ &&&& $\uparrow$ \\
$\DD$  -- Modular Fusion Category &&&&  $\CC$ -- Modular Equivariant Fusion Category\\
$\uparrow$ &&&& $\uparrow$ \\
$VOA$ && $<\mspace{-9mu}-----\mspace{-9mu}>$ && ${VOA}^{\circlearrowleft G}$ 
\end{tabularx}
\caption{Overview of Development}\label{f:development}
\end{scriptsize}
\end{figure}

The features of a usual conformal field theory were captured by the formalism of a tensor or fusion category. The notion of a Verlinde algebra which arises from a fusion category, the action of a modular group on the algebra, and the Verlinde formula gave us beautiful means to study these theories. Such they were that in \cite{Ki} Kirillov put forth  the notion of a $G$-equivariant fusion category so that it generalized the formalism of a fusion category. Roughly, a $G$-equivariant fusion category $\CC$ additionally possesses a $G$-grading and -action satisfying similar functorality properties which respect the group grading and action. Unlike the fusion category there are not braiding isomorphisms but rather, as Turaev \cite{T1} called them, $G$-crossed braidings: Denote the action of $g \in G$ on an object $V \in \CC_h$ by ${}^{g}V$. Then for $V \in \CC_g$ and $W \in \CC_h$ the isomorphism corresponding to the crossed braiding is $\sigma : V \ttt W \isoto {}^{g}W \ttt V$. Then Kirillov generalized the Verlinde algebra, giving us the so-called extended Verlinde algebra, which as he demonstrated can, like the Verlinde algebra, too have the action of the a modular group defined on it. \firef{f:development} summarizes this development.

The paper is presented in two parts. In \paref{p:first} we synthesize known results from several papers. The features of the theory needed later are highlighted and when necessary discussed at greater length. In the \paref{p:second} we continue to generalize the theory associated to modular fusion categories and their Verlinde algebras. Here we have a major break from the previous theory and, more generally, from the modules of classical objects: [as a result of the group action] the product is not commutative. In this paper we present the following main results:

\begin{itemize}
  \item We show for any modular $G$-equivariant fusion category that the fusion rules are diagonalizable for a certain subalgebra of the extended Verlinde algebra. Specifically, one can consider the space $\Ve_{\e}$ as a $\Ve_{\e,\e}$ module. In this case we diagonalize the fusion rules. We then generalize the classical Verlinde formula to the products of this type.
  \item We show for any modular $\Z_{2}$-equivariant fusion category that the fusion rules are diagonalizable for the entire algebra. Then we generalize the Verlinde formula to products in $\Ve_{\e}$.
  \item We consider an important example of a modular $\Z_{2}$-equivariant fusion category. This example arises when one considers the type $D_{2m+2}$ quantum subgroup of the semi-simple representations of $\U_{q}(\mathfrak{sl}_2)$. We give this example a thorough treatment, applying the results presented in this paper.
\end{itemize}

In \seref{s:setup} we give the definition of an equivariant fusion category. \seref{s:eqfuexamples} has trivial examples of such categories as they were offered in \cite{Ki}. \seref{s:exva} gives the definition of our main object of study, the extended Verlinde algebra. Here we provide some commentary to clarify possible points of confusion as well as a second, equivalent, definition of the extended Verlinde algebra. In \seref{s:algorb} we discuss two important constructions: the {\em algebra in a category} and the {\em orbifold construction}. The algebra in a category construction was used by Kirillov and Ostrik in \cite{KO} when they formalized the notion of a {\it quantum subgroup} and the $q$-analog of McKay Correspondence in $\U_q(\mathfrak{sl}_{2})$. The orbifold construction is presented in \cite{Ki}. These constructions provide a way for one to pass between modular fusion categories and equivariant fusion categories. The two constructions for the first time are presented along-side one another to show how they relate. These constructions will play an important role in \seref{s:uqsl} where we discuss at length the equivariant fusion category associated to the quantum subgroup of type $D_{2m+2}$. In \paref{p:second} of the paper we present new work. \seref{s:basis} deals with the lack of canonical basis in the algebra. Previously, in the Verlinde algebra, one could simply take a representative from each isomorphism class; the basis was then given by the identity map of each of these representatives. The group action ensure that there is no longer such a canonical basis and also introduces the non-commutativity in the product. In \seref{s:v11} we discuss the fusion with elements of the subalgebra [a standard Verlinde algebra]. Here we show that the fusion rules can be diagonalized and give the first generalization of the Verlinde formula. In \seref{s:uqsl} we develop an important real-life example of an equivariant fusion category. All the theory introduced in \paref{p:first} is employed in the example. We related the $s$-matrix of the two algebras via the algebra in a category construction from \cite{KO} and show explicitly that the main result from \seref{s:v11} does indeed predict the fusion rules. \seref{s:Z2}, the final section of the paper, considers the case $G = \Z_{2}$. Here we show that one can indeed diagonalize the fusion rules. We give another generalized version of the Verlinde formula and we consider the example from \seref{s:uqsl} in this context.

\begin{figure}[h]
  \vspace{20pt}
  \raisebox{-0.5\height}{\begin{overpic}
  {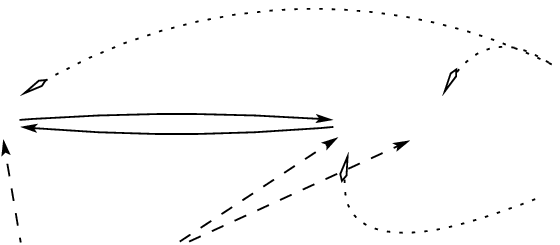}
  \put(-3,20){$\DD$}
  \put(63,20){$\CC \supset \CC_{\e}$}
  \put(2,-9){$V^{G} \subset V^{\circlearrowleft^{G}}$}
  \put(20,14){\sz{orbifolding}}
  \put(13,26){\sz{algebra in category}}
  \put(-40,20){\fz{Modules}}
  \put(-40,-9){\fz{VOAs}}
  \put(90,28){\sz{Fusion Categories}}
  \put(82,10){\sz{Equivariant Fusion Category}}
  \end{overpic}}
  \vspace{13pt}
  \caption{Fusion and Equivariant Fusion Categories}\label{f:roadmap1}
\end{figure}

It is known that for most cases one cannot achieve the same results. That is, even for most abelian groups the fusion rules cannot be diagonalized. In a future paper I will extend these results to cyclic groups. I am beginning to classify these categories and I will include more general results in that paper. As another important example I will also consider modules over twisted affine Lie algebras.

\part{General Setup}\label{p:first}
\section{Equivariant Fusion Categories}\label{s:setup}

The definitions and results in this section were already established by Kirillov in \cite{Ki}. To define the extended Verlinde algebra we will need a $G$-equivariant fusion category. The definition is due to Turaev. See \cite{T1}. As formulated in \cite{Ki} let us recall this

\setcounter{subsection}{-1}

\begin{definition}\label{d:gequ}
  A {\sem $G$-equivariant category $\CC$} is an abelian category with the
following additional structure:

\begin{description}
\item[$G$-grading] Decomposition of $\CC$ over the group $G$.
\[
\CC=\bigoplus_{h \in G} \CC_h
\]

where each $\CC_h$ is a full subcategory in $\CC$. 
\\
\item[Action of $G$] For each $g\in G$, we are given a  functor
  $\Rho_g\colon \CC\to \CC$ and functorial isomorphisms
  $\alpha_{gh}\colon \Rho_g\circ \Rho_h\isoto \Rho_{gh}$ such that
  $\Rho_\e=\id$, $\Rho_g \CC_h\subset \CC_{ghg^{-1}}$, and $\alpha_{g_1g_2,
      g_3}\circ \alpha_{g_1,g_2} =\alpha_{g_1, g_2g_3}\circ \alpha_{g_2, g_3}$
  [both sides are functorial isomorphisms 
  $\Rho_{g_1}\Rho_{g_2}\Rho_{g_3}\isoto \Rho_{g_1g_2g_3}$ ].

\end{description} 

We will use the notation ${}^gV$ for $\Rho_g(V)$.
\end{definition}

\begin{definition}\label{d:geqfu}
A {\sem $G$-equivariant fusion category} is a semisimple $G$-equivariant
abelian category with additional structure. As follows:

\begin{itemize}
\item The structure of a rigid monoidal category such that\\
  $\1$ is a simple object\\
  $\Rho_g$ is a tensor functor\\
   for $X \in \CC_g, Y \in \CC_h, X \otimes Y \in \CC_{gh}$
   \\
\item Functorial isomorphisms $\delta_{V}\colon V\to V^{**}$, satisfying
  the same compatibility conditions as in the absence of $G$ (see
  \cite{BK}) and the additional condition
  $\Rho_g(\delta_V)=\delta_{\Rho_g(V)}$. 
  \\
\item A collection of functorial isomorphisms $R_{V,W}\colon V\ttt
  W\to {}^gW\ttt V$ for every $V\in \CC_g, W\in \CC_h$, satisfying
  an analog of the pentagon axiom (see \cite[Section 2.2]{T1}). 
\end{itemize}
\end{definition}

In our study $\CC$ shall always refer to a modular $G$-equivariant fusion category where $G$ is a finite group. $\DD$ shall always refer to a fusion category. Other categories will be denoted by other symbols.

The category has the following properties which are immediate from the definition. 

\begin{proposition} Let $\CC$ be a $G$-equivariant fusion category. Then
\begin{enumerate}
  \item The unit object $\1 \in \CC_{\e}$.
  \item For $V \in \CC_g$ we have that $V^* \in \CC_{g^{-1}}$.
  \item By rigidity, we have ${}^g1 = 1$
  \item and $({}^gV)^* = {}^g(V^{*})$ canonically.
\end{enumerate}
\end{proposition}

The $G$-equivariant fusion category $\CC$ is not a fusion category. The braiding in a fusion category is commutative while the braiding in $\CC$ is not. Note that a fusion category is a special case when the grading of the category is trivial. The full subcategory $\CC_{\e} \subset \CC$, commonly known as the {\it untwisted sector}, is a fusion category. Recall \firef{f:roadmap1}. The category $\CC_{\e}$ will have an important role throughout this paper.

As in the case of a fusion category the existence of functorial isomorphisms $V \simeq V^{**}$ allows us to define a system of twists $\theta_V$ in the category.

\begin{lemma}\label{l:twists} Let $\CC$ be a $G$-equivariant fusion category. Then $\CC$ has a collection of functorial isomorphisms $\theta_V \colon V \to {}^gV$ for $V \in \CC$ with the following properties
  \begin{enumerate}
    \item $\theta_{\1} = id$ 
    \item $\theta_{U \ttt V} = (\theta \ttt \theta) \Rho_{{}^gV,U} \Rho_{U,V}$
    \item $\theta_{V^{*}} = \Rho_{g^{-1}} (\theta^*_V)$
    \item $\theta_{{}^hV} = \Rho_h (\theta_V)$
  \end{enumerate}
\end{lemma}
  
There is a graphical calculus which has been developed to represent morphisms in $\CC$. This generalizes the technique of representing morphisms in a braided tensor category by tangles. This is the work of Turaev \cite{T1}. The uninitiated reader may refer to \cite{Ki} for a brief review of this technique.

\section{Examples of Equivariant Fusion Categories}\label{s:eqfuexamples}

Here we give two examples of equivariant fusion categories. Although both are trivial in some sense they do provide a nice initiation to the theory of these categories. I duplicate these examples as they were given in \cite[Section 5]{Ki}. Some details of the orbifold construction and related theorems will be given in \seref{s:algorb}. For a fuller treatment see \cite{Ki}.

\subsection{Graded Vector Space} Let $\CC$ be the category of $G$-graded vector spaces. This is the category with simple objects $X_{g}$ for $g \in G$. The tensor product is given by $X_g \ttt X_h = X_{gh}$ and duality by $X_{g}^{*} = X_{g^{-1}}$. The action of $G$ is defined by $R_{g}X_{h} = X_{ghg^{-1}}$. Then the orbifold category [see \seref{s:algorb} for definition] $\CC / G$ is the category of finite-dimensional modules over the Drinfeld double $D(G)$ with the standard tensor product. 

\subsection{Twisted Graded Vector Space} Let $\CC^{\omega}$ be a twisted category of $G$-graded vector spaces. Here the category is a rigid monoidal category that coincides with the category given above. Tensor product and duality are instead defined as follows: $X_{g} \ttt X_{h} \isoto X_{gh}$, and $X_{g}^{*} \isoto X_{g^{-1}}$ [non-canonically]. This category defines a $3$-cocycle $\omega \in \CC^{3}(G,\C^{\times})$. Suppose that we choose a system of isomorphisms $\alpha_{gh} \colon X_{g} \ttt X_{h} \to X_{gh}$. Then $\alpha_{g_1,g_2g_3} \alpha_{g_2,g_3} = \omega(g_1,g_2,g_3) \alpha_{g_1g_2,g_3} \alpha_{g_1,g_2}$. We note that two such categories are equivalent as monoidal categories if and only if $[\omega] = [\omega']$. This defines a bijection between equivalence classes of twisted categories of $G$-graded vector spaces and $H^{3}(G,\C^{\times})$. Now the action of the group can be defined by $R_{g}X = X_{g} \ttt X \ttt X_{g}^{*}$, and the braiding isomorphism as composition

	\[
	X_{g} \ttt X_{h} \isoto X_{g} \ttt X_{h} \ttt X_{g}^{*} \ttt X_{g} = {}^{g}X_{h} \ttt X_{g}.
	\]

	This defines on our category $\CC^{\omega}$ the structure of a $G$-equivariant fusion category. The corresponding orbifold category $\CC^{\omega} / G$ coincides with the category of modules over the twisted Drinfeld double $D^{\omega}(G)$. For the definition of these modules see \cite{DPR}, and \cite{DN}.

\section{From Fusion To G-Equivariant Fusion And Back}\label{s:algorb}

Before we study the extended Verlinde algebra we introduce some important categorical constructions that we will later need. The constructions and theorems presented in this section are not new. Rather, what we are doing here is synthesizing results across several papers.

First we recall the definition of a $\DD$-algebra $A$. Then we define the category of modules $\rep A$ of this algebra. See \cite{KO} for details.

\begin{definition}
  Let $\DD$ be a fusion category. An {\dem associative commutative algebra $A$ in $\DD$} is an object $A \in \DD$ along with morphisms $\mu \colon  A \ttt A \to A$ and $\iota_{A} \colon  \1 \injto A$ satisfying associativity and commutativity compatibility conditions along with a unique unit element.

\end{definition}

We define the category of modules of this algebra as follows.

\begin{definition}\label{d:repa}
      Let $\DD$ be a fusion category and $A$ a $\DD$-algebra. Define the category $\rep A$ as follows: 
      \begin{itemize}
	\item The objects are pairs $(V, \mu_V)$ where $V\in \DD$ and $\mu_V\colon A\ttt V\to V$ is a morphism in $\DD$ satisfying the following properties:
      \begin{enumerate}
      \item $\mu_V \circ (\mu\ttt \id)=\mu_V \circ (\id \ttt \mu_V)\colon A\ttt A\ttt V\to V$
      \item $\mu_V(\iota_A \ttt \id)=\id\colon \1 \ttt V\to V$
      \end{enumerate}
	\item The morphisms are defined by 
  \begin{multline*}
  \Hom_{\rep A}((V, \mu_V), (W, \mu_W))\\
  =\{ \ph\in \Hom_\C(V, W) \vert \mu_W\circ(\id \ttt \ph)=\ph\circ\mu_V\colon A\ttt V\to W\}
  \end{multline*}
  \end{itemize}
\end{definition}

We say that an algebra $A$ in $\DD$ is {\dem rigid} if $\1$ has multiplicity $1$ in $A$ and the composition $A \ttt A \xxto{\mu} A \to \1$ is a non-degenerate pairing. See \cite{KO} for details.

\begin{figure}[th]
 \vspace{10pt}
  \raisebox{-0.5\height}{\begin{overpic}
  {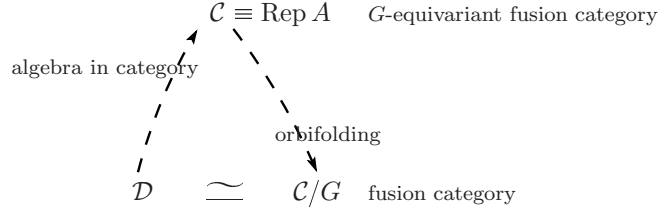}
  \put(-2,-14){$\DD$}
  \put(40,83){$\CC \equiv \rep A$}
  \put(85,-14){$\CC /G$}
  \put(-67,54){\fz{algebra in category}}
  \put(75,18){\fz{orbifolding}}
  \put(37,-17){$\isomor$}
  \put(125,83){\fz{$G$-equivariant fusion category}}
  \put(125,-14){\fz{fusion category}}
  \end{overpic}}
  \vspace{10pt}
\caption{Algebra in a Category and Orbifolding}\label{f:roadmap2}
\end{figure}

\begin{theorem}\label{t:algcat} Let $\DD$ be a fusion category. Suppose that $A \in \DD$ is an associative commutative algebra with the following properties:
  \begin{itemize}
    \item $\theta_{A} = \id$
    \item $A$ is rigid
    \item There is an action of a finite group $G$ by automorphisms $\pi_{g}$ of $A$ such that the action is faithful and $A^{G} = 1$.
  \end{itemize}
  Then the category $\CC \equiv \rep A$ is a $G$-equivariant fusion category. The objects and morphisms are given in \deref{d:repa}.
\end{theorem}

See \cite[Section 4]{Ki} and \cite[Section 5]{orbi2} for the details.

Between the categories $\DD$ and $\CC$ define two functors $F \colon \DD \to \CC$ and $G \colon \CC \to \DD$ as follows
  \[
  \begin{split}
  F & (V) = A \ttt V \mspace{3mu}, \mspace{6mu} \mu_{F(V)} = \mu \ttt \id \\
  \text{and} \mspace{9mu} G & ((V,\mu_{V})) = V.
  \end{split}
  \]

We have the following important theorems from \cite{KO}

\begin{theorem}
$F$ and $G$ are exact and injective on morphisms. $F$ and $G$ are adjoint. $F$ is a tensor functor, $F(V \ttt W) = F(V) \ttt_{A} F(W)$ and $F(\1) = A$. And $G(F(V)) = A \ttt V$.
\end{theorem}

\begin{theorem}\label{t:dimrel}
  For $X,Y \in \CC \equiv \rep A$ the following hold 
  \begin{align*}
  dim_{\CC}(X) &= \frac{dim_{\DD}(X)}{dim_{\DD}A},\\
  \text{\rm{and }} dim_{\CC}(F(V)) &= dim_{\DD}(V).
  \end{align*}
\end{theorem}

Let us recall the orbifold construction. A detailed account is given in \cite{orbi2}. 

\begin{definition}
  Let $\CC$ be a $G$-equivariant fusion category. Then the orbifold category $\CC / G$ is defined as follows:
  \begin{itemize}
    \item The objects are pairs $(X,\Phi)$ where $X \in \CC$ and $\Phi = \{ \ph_{g} \}_{g \in G}$ is a collection of $\CC$-morphisms $\ph \colon {}^{g}X \isoto X$ such that $\ph_{\e} = \id$ and $\ph_{g}R_{g}(\ph_{h}) = \ph_{gh}$.
    \item The morphisms $(X,\Phi) \to (Y,\Psi)$ are $\CC$-morphisms $\tau \colon X \to Y$ such that $\psi_{g} \circ R_{g}(\tau) = \tau \circ \ph_{g}$ for all $g$ in $G$.
  \end{itemize}
\end{definition}

\begin{nb} Not every object $X \in \DD$ will have such a set $\Phi$ of morphisms. In which case the object does not appear as part of a pair in the orbifold construction. Note also that for a given $X \in \CC$ there may be more than one such set $\Phi$.
\end{nb}

The relation between the algebra in a category construction and the orbifold construction is summarized in \firef{f:roadmap2} and given in \cite{orbi2} by the following

\begin{theorem}
Let $\DD$ be a fusion category and $A$ a commutative algebra in $\DD$ satisfying the conditions given in \thref{t:algcat}. Then the category $\DD$ is naturally equivalent to the orbifold category $\rep A / G$.
\end{theorem}

Like before, we can define adjoint functors between the two categories $\CC$ and $\CC / G$. Define $F' \colon \CC / G \to \CC$  as follows
  \[
  F' ((X,\{ \ph_{g} \})) = X.
  \]

\begin{figure}[bh]
 \vspace{10pt}
  \raisebox{-0.5\height}{\begin{overpic}
  {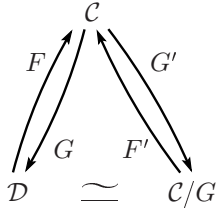}
  \put(-2,-14){$\DD$}
  \put(40,83){$\CC$}
  \put(85,-14){$\CC /G$}
  \put(8,56){$F$}
  \put(23,9){$G$}
  \put(60,9){$F'$}
  \put(75,56){$G'$}
  \put(37,-17){$\isomor$}
  \end{overpic}}
  \vspace{10pt}
\caption{Functors between the Categories}\label{f:roadmap3}
\end{figure}

Define $G' \colon \CC \to \CC / G$ as follows
  \[
  G' (V) = (X,\{ \ph_{g} \})
  \]
where $X = \oplus_{h \in G} {}^hV$ and $\ph_{g} \colon \oplus_{h} {}^{gh}V \to \oplus_{h} {}^{h}V$. These functors are summarized in \firef{f:roadmap3}.

\section{The Extended Verlinde Algebra}\label{s:exva}
\setcounter{subsection}{-1}

Let us begin by recalling the definition of the Verlinde algebra.

For $\DD$ a fusion category denote its Grothendieck ring by $K(\DD)$. The algebra $\V(\DD) = K(\DD) \ttt_{\Z} \C$ is a finite-dimensional commutative associative algebra. The algebra has a basis $V_i = \langle V_i \rangle$ , $i \in I$ and a unit $1 = \langle V_0 \rangle$. This algebra is called the Verlinde algebra. Next we give an equivalent definition of the Verlinde algebra that we will later generalize.

\begin{definition} Let $\DD$ a fusion category. The {\dem Verlinde algebra} $\V = \V(\DD)$ is the complex vector space given by
    \[
    \V = \bigoplus_{i \in I} \Mor_{\DD}(V_i , V_i) ,
    \]
    where the sum is over the set $I = I(\DD)$ of isomorphism classes of simple objects in $\DD$. The algebra $\V$ has a natural basis. For each isomorphism class fix a representative $V_i$. The canonical basis of $\V$ is the identity map $\lambda \colon V_i \to V_i$ of each of the representatives.
\end{definition}

The Verlinde algebra was generalized by Kirillov \cite{Ki}. The generalization was motivated by the modular functor approach: the generalized algebra can be defined as a vector space associated to a torus with no punctures. Kirillov gives us the the following

\begin{definition}\label{d:xverlinde}
Let $\CC$ be a $G$-equivariant fusion category. For each isomorphism class of simple objects in $\CC$ fix a representative $V_i$. Then the {\sem extended Verlinde algebra} of $\CC$ is the complex vector space defined by

	\[
	\Ve(\CC)=\bigoplus_{i\in I, g\in G}\Mor_{\CC}(V_i, {}^g V_i).
	\]
	
\end{definition}

In our study we consider extended Verlinde algebras $\Ve(\CC)$ where the category $\CC$ is such that $G$ is a finite group and there are finitely many isomorphism classes of simple objects. Thus our extended Verlinde algebra is finite dimensional.

\begin{notation} The identity element of the group is written $\e$. We index the classes of simple objects by the set $I = I(\CC)$. $I_h$ is used to denote the classes of simple objects in $\CC_h$. Lastly, ${}^gI_h$ shall denote the classes in $I_h$ which are invariant under the action of $g$. We use the notation $V_i$ to denote the representative that we have fixed. In particular $V_0$ where $0 \in I$ is the unit object $\1$ in $\CC$.
\end{notation}

Note that for $V_i \in \CC_h$ the space $\Mor_{\CC}(V_i, {}^gV_i)$ is empty unless $gh=hg$. Thus the extended Verlinde algebra can be written as

	\[
	\Ve(\CC)=\bigoplus_{g,h|gh=hg} \Ve_{g,h}(\CC); \quad 
	\Ve_{g,h}(\CC)=\bigoplus_{i\in I_h}\Mor_{\CC}(V_i, {}^g V_i).
	\]

In particular, $\Mor(V_i, {}^g V_i)$ is trivial unless ${}^g V_i$ belongs to the same isomorphism class as $V_i$. So we note that there may be $i \in I_h$ such that $\Mor(V_i,{}^gV_i)$ is trivial although the condition $gh=hg$ is satisfied. Stated differently, the space $\Mor(V_i, {}^g V_i)$ is non-trivial if and only if $V_i$ as an isomorphism class is invariant under the action of $g$. Equivalently, $\Mor(V_i,{}^gV_i)$ where $i \in I_h$ is non-trivial if and only if $i \in {}^gI_h$. We introduce another formulation of the extended Verlinde algebra:
	
	\[
	\Ve(\CC)    =\bigoplus_{g \in G} \Ve_{g}(\CC); \quad 
	\Ve_{g}(\CC)=\bigoplus_{i \in {}^{{}^g}\!I}\Mor(V_i,{}^gV_i)	
	\]

	This second formulation of the definition puts emphasis on the invariance of a class under the action of the group. In this formulation the direct sum has no trivial summands. That is, $\Mor(V_i,{}^gV_i)$ is non-trivial for all $g \in G$ and $i \in {}^gI$ since $\1 \isomor {}^g\1$ for all $g$ in $G$ and we consider elements of the set ${}^gI$ rather than $I$. The indexing set ${}^{g}I$ can be written as a disjoint union of sets: ${}^gI = \sqcup\medspace {}^gI_h$ over $h \in H_g \subset G$. The subset $H_g$ is determined by the group action on the set of isomorphism classes. From the space $\Ve_{g}$ we can in a natural way recover the spaces $\Ve_{g,h}$. 

\begin{nb} Suppose that  $X \isoto {}^gX$ in the category. Then the morphism $\ph \colon X \to {}^gX$ does not necessarily appear as an element of the algebra. For example, $X = V_i \oplus V_j$, where ${}^gV_i \backsimeq V_j$ and ${}^gV_j \backsimeq V_i$. The obstruction is that $X$ must be the direct sum of $g$-invariant [up to isomorphism] simple objects in the category. \end{nb} 

Let us introduce a formulation of the algebra that is based on isomorphism classes. Then with this formulation in place we will be prepared to discuss further the properties of the algebra.
	
\subsection{Classes and the extended Verlinde algebra}\label{s:classes}

$\Ve_g$ is isomorphic to the vector space spanned by classes $[\ph]$ where $\ph \colon V \to {}^g V$, $V$ is a linear combination of $V_i$ for $i \in {}^gI$, with the following relations 

\begin{enumerate}
  \item For any $\lambda \in \C$ $\ph \colon V \to {}^gV$ and $\psi \colon V \to {}^gV$ one has
    \[
    \lambda [ \ph ] = [ \lambda \ph ] \quad \textrm{and} \quad  [\ph] + [\psi] = [\ph + \psi].
    \]
  \item For any $\ph \colon V \to {}^gV$ and isomorphism $T \colon V \isoto V^{'}$ one has
    \[
    [R_g(T)\ph T^{-1}] = [\ph].
    \]
  \item Suppose $V = \oplus V_i$ for some $i \in {}^gI$ and $\ph \colon V \to V$ is given by $\ph = \sum \ph_{ij}$ where $\ph_{ij} \colon V_i \to {}^gV_j$. Then
    \[
    [\ph] = \sum [\ph_{ii}].
    \]
    
\end{enumerate}

We are now prepared to talk some about the structure of the extended Verlinde algebra. For a detailed treatment see \cite[Section 8]{Ki}.

\subsection{Extended Verlinde algebra structure}\label{s:struc}
\begin{description}
  \item[Dimension] The dimension $d_i$ of an object $V_i$ in $\CC$ is defined to be the categorical trace of the identity of $V_i$. 
 \\ 
  \item[Tensor Product] The extended Verlinde algebra is an associative algebra. We denote the product as $\ttt$. Suppose $\ph \in \Ve_{g,a}$ and $\psi \in \Ve_{h,b}$. Then the product when $g = h$ is defined to be
    \[
    [\ph] \ttt [\psi] = [\ph \ttt \psi] \in \Ve_{g,ab}
    \]
    and is zero otherwise. The tensor product is in $\Ve_{g,ab}$ because the action of the group is a tensor functor. When we use the word fusion or fusion product we are referring to the tensor product of the algebra.
\\
  \item[Convolution Product] The algebra has another associative product. We denote it by $*$. Let $\ph \colon V_i \to {}^gV_i$ and $\psi \colon V_j \to {}^hV_j$. The convolution product $[\ph] * [\psi]$ is zero when $\, {}^hV_j$ and $V_i$ are not isomorphic. Otherwise, let $\kappa$ be an isomorphism between ${}^hV_j$ and $V_i$. Then
    \[
    [\ph]*[\psi] = d_i^{-1} [V_j\xxto{\psi}{}^{h}V_j\xxto{\kappa}  
    V_i\xxto{\ph} {}^{g}V_i\xxto{\Rho_{g}(\kappa^{-1})} {}^{gh}V_j ]
    \]
    where $d_i$ is the dimension of $V_i$. Thus for $\ph \in \Ve_{g,a}$ and $\psi \in \Ve_{h,a}$ we have that $[\ph] * [\psi] \in \Ve_{gh,a}$. The coefficient is in the definition to make life simpler later on. Note that we have for a simple object $V$ and $\ph \colon {}^{h}V \to {}^{gh}V$ and $\psi \colon V \to {}^hV$ that $[\ph] * [\psi] = \frac{1}{dim \, V}[\ph \psi]$.
  \\  
  \item[Bilinear Form] To define the bilinear form we first make some preliminary definitions. Recall that for $\ph \colon V \to {}^gV$ we have an adjoint morphism $\ph^* \colon {}^gV^* \to V^*$. This defines on $\Ve$ a linear map $*\colon \Ve_{g,h}\to \Ve_{g^{-1}, h^{-1}}$ such that $(\ph \ttt \psi)^* = \psi^* \ttt \ph^*$ and $(\ph * \psi)^* = \psi^**\ph^*$.
  
Next, define the constant term map $[ \ ]_0 \colon \Ve \to \C$ as follows
  \begin{align*}
    &[\ph]_0=0, \quad \ph\colon V_i\to {}^gV_i, \ i\ne 0\\
    &[\chi^g_0]_0=1
  \end{align*}

  where $\chi_0^g\colon \1\to {}^g\1$ is the canonical
  isomorphism. One easily sees that it completely determines $[\ ]_0$
  and that $[x]_0=0$ if $x\in \Ve_{g,h}, h\ne 1$.

  We can now define the bilinear form 
  \[
  (\ph,\psi)=[\ph\ttt \psi^*]_0.
  \]
\end{description}

\vspace{3pt}

\begin{proposition} The bilinear form enjoys the following properties: 
    \begin{enumerate}
      \item 
	For $\ph\in \Ve_{g_1,h_1}, \psi\in \Ve_{g_2,h_2}$, we have 
	$(\ph,\psi)=0$ unless $g_1=g_2^{-1}, h_1=h_2$.
      \item The form is symmetric: $(\ph,\psi)=(\psi,\ph)$, non-degenerate,
	  and $G$-invariant. 
	\item $(\chi_i,\chi_j)=\de_{ij}$
	\item $(x\ttt y,z)=(x,z\ttt y^*)$. 
      \end{enumerate}
    \end{proposition}

In the extended Verlinde algebra the tensor product and the convolution product are not commutative. In general the products fail to be commutative even when the group is abelian. This is unlike the case of the Verlinde algebra where both products commute.

As in the case of the Verlinde algebra it is possible to define linear operators on the extended Verlinde algebra which under some non-degeneracy condition defines the action of the modular group $SL_{2}(\Z)$ on $\Ve$. We introduce these operators now.

\begin{definition} $\tilde t \colon \Ve_{g,h} \to \Ve_{gh,h}$ is defined by
  \[
  [\ph] \mapsto [\theta \ph] = [\ph \theta]
  \]
  where $\theta$ is the universal twist.
\end{definition}

The equality $[\theta \ph] = [\ph \theta]$ follows from \leref{l:twists} and \seref{s:classes}. The existence of this operator will later be used to develop theorems on the structure of particular extended Verlinde algebras. Let us now define a most useful operator on the extended Verlinde algebra.

\begin{definition}\label{d:smatrix}
  The linear operator $\tilde s \colon \Ve \to \Ve$ called the {\sem $s$-matrix} is defined as follows.
  \begin{equation*}\label{e:smatrix}
    \begin{aligned}
      \tilde s \colon \Ve_{g,h} &\to \Ve_{h^{-1},g} \\
      \tilde s[\varphi] &=\mspace{-9mu}\sum_{k \, \in \, \vphantom{I}^{h^{-1}}\!I_g}\mspace{-9mu}(\tilde s[\varphi])_k
    \end{aligned}
  \end{equation*}
  where, for $\varphi\colon V\to {}^g V, V\in \CC_h$, we define  $(\tilde s[\varphi])_k \colon V_k \to {}^{h^{-1}}V_k, V_k \in \CC_g$ by \firef{f:smatrix}.
  
  \begin{figure}[th]
    $
    (\tilde s[\varphi])_k=d_k\mspace{9mu}
    \raisebox{-0.5\height}{\begin{overpic}
      {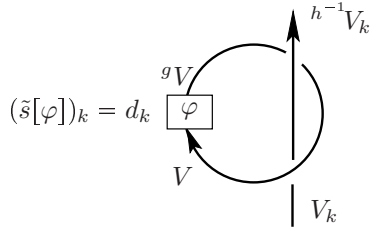}
      \put(7,52){$\varphi$}
      \put(3,20){$V$}
      \put(-2,65){${}^gV$}
      \put(65,5){$V_k$}
      \put(65,90){${}^{h^{-1}}\!V_k$}
    \end{overpic}}
    $
    \caption{$s$-matrix}\label{f:smatrix}
  \end{figure}

\end{definition}

When using the graphical calculus to represent objects in the algebra it will be useful to have the following 

\begin{lemma} Let $\ph \in \Ve_{g,h}$ and $\psi \in \Ve_{h,g}$. Suppose $\ph \colon V \to {}^gV$ and $\psi \colon W \to {}^hW$. Then

  \[
  (\tilde s\ph,\psi)=
  \raisebox{-0.5\height}{\begin{overpic}
  {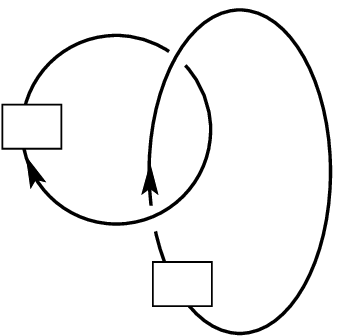}
  \put(5,60){$\ph$}
  \put(52,13){$\psi$}
  \put(10,28){$V$}
  \put(5,85){${}^gV$}
  \put(85,40){$W$}
  \put(53,25){${}^hW$}
  \end{overpic}}
  \]

\end{lemma}

The operators $\tilde s$ and $\tilde t$ behave nicely with the bilinear form.

\begin{lemma}\label{l:sprops} The operators $\tilde s$ and $\tilde t$ have the following properties.
  \begin{enumerate}
    \item When restricted to $\Ve_{1,1}$ the operators $\tilde s$ and $\tilde t$ are the corresponding operators of the [non-extended] Verlinde algebra $\V(\CC_1)$.
    \item $\tilde s$ and $\tilde t$ are symmetric: $(\tilde s \ph , \psi ) = (\ph , \tilde s \psi)$ and $(\tilde t \ph , \psi ) = (\ph , \tilde t \psi)$.
    \item $(\tilde s \ph , \psi) = (\tilde s \psi^* \! , \ph^*)$.
    \item $\tilde s ( \ph \ttt  \psi ) = \tilde s (\psi) * \tilde s (\ph)$.
  \end{enumerate}
\end{lemma}

\begin{proof}
  Immediate from the graphical calculus.
\end{proof}

The operator $\tilde s$ enjoys the property of interchanging the tensor product with the convolution product. We will leverage this fact with the structure of the extended Verlinde algebra to study the tensor product explicitly. Later still we will make use of the group action in our study of the tensor product.

\begin{definition}\label{d:modular}
  A $G$-equivariant fusion category with finitely many isomorphism classes of simple objects is called {\sem modular} if the operator $\tilde s$ is invertible.
\end{definition}

The category is called modular since, after re-normalization, the operators $\tilde s$ and $\tilde t$ satisfy the relations of $SL_2(\Z)$. See \cite[Section 10]{Ki} for details. In what follows it will be convenient to have the $s$-matrix normalized so that it is unitary.

\begin{proposition}\label{p:norm}
  Define the numbers $p^{\pm}$ as follows
  \[
  p^{\pm} = \sum_{i \in I_{\e}} \theta_{i}^{\pm 1} d_{i}^{2}.
  \]
  Put 
  \[
  D = \sqrt{p^{+} p^{-}} \mspace{3mu}, \mspace{9mu} s = D^{-1} \tilde s.
  \]

Then $s$ is a symmetric unitary operator.
\end{proposition}

The $s$-matrix of a modular fusion category can also be normalized so that it is unitary. The number $D=D^{\DD}$ is defined by the same formula as in the proposition.

\begin{nb}
The equality \( \tilde s (\ph \ttt \psi) = \tilde s(\psi) * \tilde s(\ph) \) from \leref{l:sprops} becomes \( s(\ph \ttt \psi) = D s(\psi) * s(\ph) \). 
\end{nb}

Modularity is preserved when moving between fusion categories and equivariant fusion categories via the algebra in a category construction or the orbifold construction. See \firef{f:roadmap2}. From \cite[Section 10]{Ki} we have the following

\begin{theorem}\label{t:modmod} 
  A $G$-equivariant fusion category $\CC$ is modular if and only if the orbifold category $\CC / G$ is modular.
\end{theorem}

So by the equivalence of the categories $\DD$ and $\CC / G$ in \firef{f:roadmap2} we have that modularity is preserved. Also note that if $\CC$ is modular then $\CC_{\e}$ is modular. It is important to realize that this implies that the restriction of the $s$-matrix from $\Ve(\CC)$ to the subalgebra $\Ve_{\e,\e}$ is an invertible operator.

To introduce the relation between the two operators $s^{\CC}$ and $s^{\CC / G} = s^{\DD}$ we extend the definition of $F'$ and $G'$ from the language of categories to algebras. Start by recalling the definitions of $F'$ and $G'$.
  \begin{align*}
    &F' \colon \CC / G \to \CC \\
    \text{is given by } &F' ((X,\{ \ph_{g} \})) = X ,\\
    \intertext{and}
    &G' \colon \CC \to \CC / G \\
    \text{is given by } &G' (V) = (X,\{ \ph_{g} \})
  \end{align*}

  where $X = \oplus_{h \in G} {}^hV$ and $\ph_{g} \colon \oplus_{h} {}^{gh}V \to \oplus_{h} {}^{h}V$. Note $\ph_{g} \in \Ve_{g^{-1}}$. See \firef{f:roadmap3}. The functors $F'$ and $G'$ can be extended to maps between the algebras as follows. Define 
  \begin{align*}
  &\widehat F' \colon \V(\CC / G) \to \Ve(\CC) \\
  \text{by } &\widehat F'((X,\{ \ph_{g} \})) = \sum_{g \in G} [\ph_{g}].
\end{align*}

For $\ph \colon V \to {}^{g}V$ define 
  \begin{align*}
   &\widehat G' \colon \Ve(\CC) \to \V(\CC / G) \\
   \text{by } &\widehat G'(\ph) = \oplus_{h} R_{h}(\ph) \colon {}^{h}V \to {}^{hg}V.
 \end{align*}

  The $\tilde s^{\CC / G}$-matrix of the Verlinde algebra of the orbifolded category is related to the $\tilde s^{\CC}$-matrix of the extended Verlinde algebra $\V(\CC)$ by the following equation from \cite[Section 9]{Ki}. See \firef{f:roadmap3}.

\begin{theorem}\label{t:srelat} For $x \in \V(\CC / G)$ and $y \in \Ve(\CC)$ we have
  \[
  (\tilde s \widehat F'x,y)_{\CC} = \frac{1}{\left| G \right|}(\tilde s x, \widehat G'y)_{\CC / G}.
  \]
\end{theorem}

We can normalize both operators so that they are unitary. The above relation between the normalized operators is given by the following

\begin{corollary}\label{c:srelat} For $x \in \V(\CC / G)$ and $y \in \Ve(\CC)$ we have
  \[
  (s \widehat F' x,y)_{\CC} = (s x, \widehat G'y)_{\CC / G}.
  \]
\end{corollary}

\begin{proof}
A theorem from \cite[Section 10]{Ki} gives us the relation $\left| G \right| D_{\CC} = D_{\CC / G}$. Apply the relation to \thref{t:srelat} and the result follows. 
\end{proof}

Of course, since the categories $\DD$ and $\CC / G$ are equivalent, \coref{c:srelat} can be used to relate values of $s^{\CC}$ and $s^{\DD}$. We will make tacit use of this equivalence when using this corollary in \seref{s:uqsl}.

\part{Results}\label{p:second}

\section{A Basis for the Algebra}\label{s:basis}

In the case of the standard Verlinde algebra we could choose a basis canonically. To each simple object $V_i$ in a fusion category there corresponds an identity map $\id \colon V_i \to V_i$ in the category. These maps formed a basis for the Verlinde algebra. We do not have this nicety in the extended Verlinde algebra. There is no canonical map in the space $\Mor(V_i,{}^gV_i)$ and so there is no canonical basis for the extended Verlinde algebra. This causes some difficulty when we try to study fusion in the algebra. To get around this difficulty we introduce the notion of a generic basis for the extended Verlinde algebra.

For $i$ invariant under the action of $g$ we have that $\Mor(V_i,{}^gV_i)$ is a one dimensional vector space. From each one dimensional vector space $\Mor(V_i,{}^gV_i)$ we pick an arbitrary element in that space. This element will serve as a basis for the space. However in cases where $g = \e$ we choose the identity map in $\Mor(V_i,{}^{\e}V_i)$ as our `arbitrary' element. This can be done canonically as in the case of the standard Verlinde algebra. These elements form a generic basis of $\Ve$.

\begin{notation} Recall how we set the stage. Given a $G$-equivariant category $\CC$ we first fixed a representative $V_i$ from each of the isomorphism classes. We then defined the corresponding extended Verlinde algebra $\Ve(\CC)$. Now we fix a generic basis for the algebra as described above. We make the following notation to refer to these basis elements.

\begin{description}
 \item[Basis elements in $\Ve_\e$] For the basis element that we fixed in $\Mor(V_i,V_i)$ which is the identity map we write $\lambda_i$. When we want to state that $i \in I_h$ we write $i = i_h$. The basis element is then written $\lambda_{i_h}$.
  \\
   \item[Basis elements in $\Ve_g$] For the basis element that we fixed in $\Mor(V_i,{}^gV_i)$ we write $\basis{i}$. When we want to state that $i \in I_h$ we write $i = i_h$. The basis element is then written $\basis{i_h}$.
  \end{description}
\end{notation}

In practice, the lack of canonical basis will not trouble us. When we study fusion in the algebra it is useful to have a basis $< \lambda_k >$ satisfying ${}^g\lambda_i * {}^h\lambda_i = {}^{gh}\lambda_i$. In general however there is no such basis. Even in cases where the group is abelian, as we shall see later, one cannot expect such a basis.

\section{Fusion with an element of $\Ve_{1,1}$}\label{s:v11}

 Let us recall our goal. Given two elements in $\Ve$ we want to write their product as a sum of elements. Two non-trivial elements of our algebra have a tensor product which is non-zero when both belong to the same subspace $\Ve_g$. Fix a basis for the algebra as described above. The basis elements are morphisms of the simple objects in the category. It suffices to understand fusion on the basis elements to realize the goal. The fusion of two basis elements is then to be written as a linear combination of basis elements. Again, we denote our basis of $\Ve_g$ by $\basis{i} \colon V_i \to {}^gV_i$ where $i$ runs over the indexing set ${}^{g}I$. Then what we are looking for is ${}^{g}L^{k}_{ij}$ such that

\[
\basis{i} \ttt \basis{j} = \sum_k \, {}^{g}L^{k}_{ij} \, \basis{k} .
\]

The $L^{k}_{ij}$ are what we call fusion coefficients. For a fixed subspace $\Ve_g$, a fixed basis, and a fixed basis element $\basis{i}$ we can consider the matrix $L_i = L^{k}_{ij} = {}^gL^{k}_{ij}$  where $j$ and $k$ are over ${}^{g}I$. We shall call this matrix the matrix of left multiplication by $\basis{i}$ in the given basis. Note that for $i = i_a$ and $j = j_b$ we know a priori that $L^{k_c} = 0$ when $c \neq ab$. When it is clear from the context we suppress the $g$ superscript and write instead $L^{k}_{ij}$.

We first study fusion in $\Ve_{\e}$. Below we show without making any restrictions on our category that we can get precise results in terms of the operator $\tilde s$ in cases where one of the elements in the product belong to $\Ve_{\e,\e}$. In such cases the product is commutative.

In this section then we consider the space $\Ve_{\e}$ as a $\Ve_{\e,\e}$-module. We will first look at the special case $\Ve_{\e,\e} \times \Ve_{\e,\e} \to \Ve_{\e,\e}$. This case is nothing more than the standard Verlinde algebra presented in the language of the extended Verlinde algebra. After staging these results in the broader language we will consider the general case.

\subsection{The Verlinde Algebra}\label{s:veralg}

Let us begin by reviewing fusion in the Verlinde algebra. Here we are restricting to the fusion of elements in $\Ve_{\e,\e}$. Thus the all operators here are maps $\Ve_{\e,\e} \to \Ve_{\e,\e}$. Under these restrictions we are in the setting of a modular fusion category and the Verlinde algebra associated to it.

\begin{lemma}\label{l:verl} Let $\CC$ be a modular $G$-equivariant fusion category. Fix a basis for $\Ve(\CC)$. Fix $i \in I_{\e}$. Let $L_i$ be the operator of left multiplication by $\lambda_i$ in our basis. For $k \in I_{\e}$ define the linear operator $D_i \colon \Ve_{\e,\e} \to \Ve_{\e,\e}$ by $D_i \lambda_k = (s \lambda_i , \lambda_k) / (s \lambda_0 , \lambda_k) \, \lambda_k$. Then we have that
\[
s L_i = D_i s.
\]
\end{lemma}

\begin{proof}
    Let $\lambda_j \in \Ve_{\e,\e}$. Then from the left hand side we have
  \[
  \begin{split}
    s L_i \lambda_j &= s ( \lambda_i \ttt \lambda_j) \\
   &= D s \lambda_j * s \lambda_i
 \end{split}
  \]
  by \leref{l:sprops}. Recall the \deref{d:smatrix} of $\tilde s$. Thus 
  \[
  \begin{split}
  s \lambda_i = D^{-1} \tilde s \lambda_i &= D^{-1} \sum_{k \in I_{\e}} \, d_i \, (\tilde s \lambda_i)_k \\
  &= D^{-1} \sum_{k \in I_{\e}} \, (\tilde s \lambda_i, \lambda_k) \, \lambda_k.
  \end{split}
  \]
  Mutatis mutandis, we have the same for $\tilde s \lambda_j$. Note that $(\tilde s \lambda_0 , \lambda_i) = d_i = dim \, V_i$. Therefore
  \[
  \begin{split}
   D s \lambda_j * s \lambda_i &= D^{-1} \sum_{k \in I_{\e}} \, (\tilde s \lambda_j, \lambda_k) \, \lambda_k \, * \, \sum_{k \in I_{\e}} \, (\tilde s \lambda_i, \lambda_k) \, \lambda_k \\
   &= D^{-1} \sum_{k \in I_{\e}} \, (\tilde s \lambda_j, \lambda_k) \, \lambda_k \, * \, (\tilde s \lambda_i, \lambda_k) \, \lambda_k \\
   &= D^{-1} \sum_{k \in I_{\e}} \, \frac{(\tilde s \lambda_i, \lambda_k)}{(\tilde s \lambda_0, \lambda_k)} \, (\tilde s \lambda_j, \lambda_k) \, \lambda_k. \\
   &= D_i s \lambda_j.
 \end{split}
   \]   

   Since the convolution product is zero when $V_i$ and $V_j$ are not in the same isomorphism class we have equality between the first and second line. $\lambda_k$ is the identity map. So the convolution product is simply the composition of functions multiplied by the constant ${d_k}^{-1}$. This concludes the proof.
 
 \end{proof}

As a direct result of the previous lemma we have the well known Verlinde formula. 
 
 \begin{theorem}\label{th:verform}Verlinde Formula. Keep the hypotheses from \leref{l:verl}. Then the fusion coefficients $L_{ij}^k$ are given by the following formula.
   \[
   L_{ij}^k= \sum_p \, \frac{(s \lambda_i , \lambda_p)(s \lambda_j , \lambda_p)(s \lambda_{k}^* , \lambda_p)}{(s \lambda_0 , \lambda_p)} \mspace{3mu}, \mspace{15mu} i,j,k \in I_{\e}.
   \]
 \end{theorem}
 \begin{proof}Fix $\lambda_i$ and $\lambda_j \in \Ve_{\e,\e}$. The equation $s L_i = D_i s$ from above can be written as follows 
   \[
   \sum_r \, L_{ij}^r (s \lambda_r , \lambda_p) = \frac{(s \lambda_i , \lambda_p)(s \lambda_j , \lambda_p)}{(s \lambda_0 , \lambda_p)}.
   \]
   We multiply both sides by $(s \lambda_p , \lambda_k^*)$ and sum over $p$. Since $s$ is a symmetric unitary operator (\leref{l:sprops}) we have that $(s \lambda_r , \lambda_p)(s \lambda_p , \lambda_k^*) = (s \lambda_r , \lambda_p)(s \lambda_p^* , \lambda_{k}) = \delta_{rk}$ and the proof.
   
 \end{proof}

\begin{figure}[ht]
 \[
    d_k\mspace{9mu}
    \raisebox{-0.5\height}{\begin{overpic}
      {smatrix}
      \put(7,49){$\lambda$}
      \put(3,20){$V_i$}
      \put(-2,65){$V_i$}
      \put(65,5){$V_k$}
      \put(65,90){$V_k$}
    \end{overpic}} = (\tilde s \lambda_i , \lambda_k) \raisebox{-0.5\height}{\begin{overpic}
    {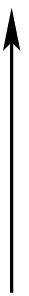}
    \put(10,5){$\lambda_{k}$}
  \end{overpic}} 
 \]
 \caption{Scalar Multiple of Identity}\label{f:multid}
\end{figure}

In the proof of \leref{l:verl} we were able to write both $(\tilde s \lambda_i)_k$ and $(\tilde s \lambda_j)_k$ as scalar multiples of the identity map $\lambda_k$. Precisely, $(\tilde s \lambda_i)_k = (\tilde s \lambda_i , \lambda_k) \, \lambda_k$. See \firef{f:multid}. The convolution product of these maps was then the composition of identity maps multiplied by a scalar. In the more general setting of $\Ve_{\e,\e} \times \Ve_{\e} \to \Ve_{\e}$ this is no longer the case. We can however still calculate the fusion coefficients in terms of the $s$-matrix. We have to adjust our procedure only slightly.

\begin{figure}[hb]
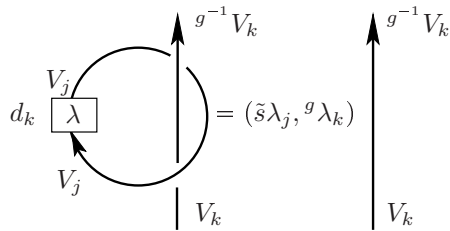

  \[
 d_k \mspace{9mu}
      \raisebox{-0.5\height}{\begin{overpic}
      {smatrix}
      \put(7,49){$\lambda$}
      \put(3,20){$V_j$}
      \put(-2,65){$V_j$}
      \put(65,5){$V_k$}
      \put(65,90){${}^{g^{-1}}V_k$}
    \end{overpic}} = (\tilde s \lambda_j, \basis{k}) \mspace{6mu}\raisebox{-0.5\height}{\begin{overpic}
    {arrow}
    \put(10,5){$V_{k}$}
    \put(10,90){${}^{g^{-1}}V_{k}$}
    \end{overpic}}
 \]
 \caption{Scalar Multiple of Basis Element}\label{f:multbasis}
 \end{figure}
 
\subsection{The Space $\Ve_{\e}$ as a $\Ve_{\e,\e}$-module}

In this section we consider $\lambda_i \ttt \lambda_j$ for $i \in I_{\e}$ and $j \in I_{g}$. Recall that the operator $\tilde s$ sends and element of $\Ve_{\e,g}$ to $\Ve_{{g^{-1}},\e}$. Thus $(\tilde s \lambda_j)_k$ is no longer a multiple of the identity map $\id \colon V_k \to V_k$. Rather $(\tilde s \lambda_j)_k$ is a multiple of ${}^{g^{-1}}\lambda_k$. See \firef{f:multbasis}. Though convolution is no longer the composition of identity maps it remains a composition of maps, one of which is an identity map. See \seref{s:struc} for the definition of the convolution product.

\begin{lemma}\label{l:v11mod}
  Let $\CC$ be a modular $G$-equivariant fusion category. Fix a basis for $\Ve(\CC)$. Fix $i \in I_{\e}$. Let $L_i \colon \Ve_{\e,g} \to \Ve_{\e,g}$ be the operator of left multiplication by $\lambda_i$ in our basis. For $k \in {}^{g^{-1}}I_{\e}$ define the linear operator $D_i \colon \Ve_{g^{-1},\e} \to \Ve_{g^{-1},\e}$ by $D_i \mspace{6mu} \vphantom{\lambda}^{g^{-1}}\mspace{-3mu}\lambda_k = (s \lambda_i , \lambda_k) / (s \lambda_0 , \lambda_k) \, \vphantom{\lambda}^{g^{-1}}\mspace{-3mu}\lambda_k$. Then we have that
  \[
  s L_i = D_i s.
  \]
\end{lemma}

\begin{proof}

Since $s$ and $D_i$ are invertible the lemma implies that $L_i$ is an invertible linear operator. We use the graphical calculus to prove the statement.

Let $\lambda_j \in \Ve_{\e,g}$. Then from the left hand side we have $s L_i \lambda_j = s ( \lambda_i \ttt \lambda_j) = D s \lambda_j * s \lambda_i$
  \[
  \begin{split}
    \\
    =&  D^{-1} \mspace{-9mu} \sum_{k \, \in \, \vphantom{I}^{g^{-1}}\!I_{\e}} \mspace{-9mu} (\tilde s \lambda_j)_k \mspace{18mu} * \mspace{18mu} \sum_{k \in I_{\e}} \, (\tilde s \lambda_i)_k \\
    =& D^{-1} \mspace{-9mu} \sum_{k \, \in \, \vphantom{I}^{g^{-1}}\!I_{\e}} \mspace{-9mu} d_k \mspace{9mu}
      \raisebox{-0.5\height}{\begin{overpic}
      {smatrix}
      \put(7,49){$\lambda$}
      \put(3,20){$V_j$}
      \put(-2,65){$V_j$}
      \put(65,5){$V_k$}
      \put(65,90){${}^{g^{-1}}V_k$}
    \end{overpic}} \mspace{27mu} * \mspace{18mu} 
    \sum_{k \in I_{\e}} \, (\tilde s \lambda_i, \lambda_k) \, \raisebox{-0.5\height}{\begin{overpic}
    {arrow}
    \put(10,5){$V_{k}$}
    \put(10,90){$V_{k}$}
  \end{overpic}} \\
   =& D^{-1} \mspace{-9mu} \sum_{k \, \in \, \vphantom{I}^{g^{-1}}\!I_{\e}} \mspace{-9mu} d_k \mspace{9mu}
      \raisebox{-0.5\height}{\begin{overpic}
      {smatrix}
      \put(7,49){$\lambda$}
      \put(3,20){$V_j$}
      \put(-2,65){$V_j$}
      \put(65,5){$V_k$}
      \put(65,90){${}^{g^{-1}}V_k$}
    \end{overpic}} \mspace{27mu} * \mspace{18mu} 
    (\tilde s \lambda_i, \lambda_k) \, \raisebox{-0.5\height}{\begin{overpic}
    {arrow}
    \put(10,5){$V_{k}$}
    \put(10,90){$V_{k}$}
  \end{overpic}} \\
    =& D^{-1} \mspace{-9mu} \sum_{k \, \in \, \vphantom{I}^{g^{-1}}\!I_{\e}} \mspace{-9mu} \frac{(\tilde s \lambda_i, \lambda_k)}{(\tilde s \lambda_0, \lambda_k)}
    \mspace{9mu} d_k \mspace{9mu}
      \raisebox{-0.5\height}{\begin{overpic}
      {smatrix}
      \put(7,49){$\lambda$}
      \put(3,20){$V_j$}
      \put(-2,65){$V_j$}
      \put(65,5){$V_k$}
      \put(65,90){${}^{g^{-1}}V_k$}
    \end{overpic}} \mspace{27mu} \\
    =& D^{-1} \mspace{-9mu} \sum_{k \, \in \, \vphantom{I}^{g^{-1}}\!I_{\e}} \mspace{-9mu} \frac{(\tilde s \lambda_i, \lambda_k)}{(\tilde s \lambda_0, \lambda_k)}
    \mspace{9mu} (\tilde s \lambda_j)_k  = D^{-1} \mspace{-9mu} \sum_{k \, \in \, \vphantom{I}^{g^{-1}}\!I_{\e}} \mspace{-9mu} \frac{(s \lambda_i, \lambda_k)}{(s \lambda_0, \lambda_k)}
    \mspace{9mu} (\tilde s \lambda_j)_k\\ 
    =& \mspace{9mu} D_i s \lambda_j.
 \end{split}
   \]   
\end{proof}

\begin{corollary}\label{c:v11mod}
  Keep the hypotheses from \leref{l:v11mod}. Then the fusion coefficients are given by the following formula. 
  \[
   L_{ij}^k = \sum_{p} \frac{(s \lambda_i, \lambda_p)(s \lambda_j, \basis{p})(s \lambda_k^*, \basis{p})}{(s \lambda_0, \lambda_p)} \mspace{3mu}, \mspace{15mu} i \in I_{\e} \mspace{6mu} \textrm{and} \mspace{6mu} j,k \in I_{g}.
  \]
\end{corollary}

Note that $p$ in sum above is indexed by the set ${}^{g^{-1}}I_{\e}$. Often we instead write ${}^gI_{\e}$ since these sets are the identical. Recall that $\Ve_{\e,g}$ is a $\Ve_{\e,\e}$ module. Thus for $\lambda_i \in \Ve_{\e,\e}$ and $\lambda_j \in \Ve_{\e,g}$ we have that $L_i \lambda_j$ is $\sum_{k} L_{ij}^k \lambda_k$ where $k \in I_g$.

\begin{proof}
Fix $\lambda_i \in \Ve_{\e,\e}$. Fix $\lambda_j \in \Ve_{\e,g}$. Then the equation $s L_i = D_i s$ can be written as follows.

\begin{align*}
    &\mspace{54mu} & s L_i \lambda_j &= D_i s \lambda_j&\\
    \iff&\mspace{54mu} & s \sum_{r \in I_g} L_{ij}^r \lambda_r &= D_i \sum_{p \in {}^{g^{-1}}I_1} (s \lambda_j, \basis{p}) {}^{g^{-1}}\lambda_p&\\
    \iff&\mspace{54mu} &\sum_{r} L_{ij}^r s \lambda_r &= \sum_{p} D_i (s \lambda_j, \basis{p}) {}^{g^{-1}}\lambda_p&\\
    \iff&\mspace{54mu} &\sum_{r} L_{ij}^r \sum_{p} (s \lambda_r, \basis{p}) {}^{g^{-1}}\lambda_p &= \sum_{p} \frac{(s \lambda_i, \lambda_p)(s \lambda_j, \basis{p})}{(s \lambda_0, \lambda_p)} {}^{g^{-1}}\lambda_p&\\
    \iff&\mspace{54mu} &\sum_{r} L_{ij}^r (s \lambda_r, \basis{p}) &= \frac{(s \lambda_i, \lambda_p)(s \lambda_j, \basis{p})}{(s \lambda_0, \lambda_p)}.& 
\end{align*}

We multiply both sides by $(s \basis{p}, \lambda_k^*)$ and sum over $p$. Since $s$ is a symmetric unitary operator we have that $(s \lambda_r, \basis{p})(s \basis{p}, \lambda_k^*) = (s \lambda_r, \basis{p})(s \basis{p}^*, \lambda_k) = \delta_{rk}$ and the proof.

\end{proof}

\subsection{Remarks} In this section we repeatedly used the fact that the structure of $\ttt$ and $*$ are isomorphic on $\Ve$. That is, that $\tilde s$ is invertible, and both $\tilde s (\ph \ttt \psi) = \tilde s (\psi) * \tilde s (\ph)$ and $\tilde s (\ph * \psi) = \tilde s (\psi) \ttt \tilde s (\ph)$ hold. Special properties that the group may impose on the algebra have yet to been used. In the proof the convolution product was $d_k^{-1} V_k \xxto{\lambda_k} V_k \xxto{{}^{g^{-1}}\lambda_k} {}^{g^{-1}}V_k$. In more general cases we see the definition of the convolution product in its full generality. Recall that $[\ph]*[\psi] = d_i^{-1} [V_j\xxto{\psi}{}^{h}V_j\xxto{\kappa} V_i\xxto{\ph} {}^{g}V_i\xxto{\Rho_{g}(\kappa^{-1})} {}^{gh}V_j ]$. In general using the $s$-matrix to interchange the tensor product with the convolution product will not ``diagonalize'' the fusion rules. Further description of $\Ve$ requires that we make use of properties imposed by the group. In cases where the group $G$ is commutative we still manage to get nice results.

Before we consider more general cases let us rest on an example. In the next section we introduce most of the machinery necessary to discuss the example. After that we consider a non-trivial example. We use the theorems from this section to compute the fusion rules. Perhaps this will lend insight into more general cases. Succinctly, a nice example is nice.

\section{An Example: Algebra of Type $D_{2m+2}$}\label{s:uqsl}

The subspace $\Ve_{\e,\e}$ is a $\Ve_{\e}$ module. In \seref{s:v11} we diagonalized the fusion rules for $\ph \ttt \psi$ in the case that $\ph \in \Ve_{\e,\e}$ and $\psi \in \Ve_{\e}$. In this section we give a non-trivial example of a modular $G$-equivariant fusion category. We start with a modular fusion category $\DD$ and consider a modular $G$-equivariant fusion category $\CC$ that arises from an algebra in $\DD$. Recall \firef{f:roadmap2}. This example is attractive to us for several reasons. For now, we show that \coref{c:v11mod} correctly predicts the fusion rules and we show how the $s$-matrix of the algebras $\V(\DD)$, $\Ve(\rep A)$ bring together the Verlinde formula with the formula [extended Verlinde formula] from our corollary.

Let $\DD$ be the semisimple part of the category of representations of $\U_q(\mathfrak{sl}_2)$ with $q = e^{\pi i / \varkappa}$ and $\varkappa \geqq 2$. This category is well understood; for a review see \cite{BK}. Denote the irreducible representations by $V_0, \dots, V_{\delta}$, where $\delta = \varkappa - 2$. These representations are the simple objects in the category. Denote the Verlinde algebra of this category by $\V = \V(\DD)$. Denote the $s$-matrix operator associated to this algebra by ${\tilde s}^{\DD}$. This category is a modular fusion category. As discussed above this algebra has a natural basis. Denote the identity map of $V_{i}$ by $\chi_{i}$. Then the $\chi_{i}$ form a basis for the algebra, where $i$ is indexed by the set $I = I(\DD) = \{0,1, \ldots, \delta \}$. Fusion in the algebra is given by the following theorem.

\begin{theorem}\label{t:quantcg}
For $\chi_{i}, \chi_{j} \in \V$. Then
  \[
  \chi_{i} \ttt \chi_{j} = \sum_{k} N^{k}_{ij} \chi_{k}
  \]
where
  \[
  N^{k}_{ij} =\begin{cases}
    1&  \text{ for $\lvert i-j \rvert \leqq k \leqq i+j$, $k \leqq 2\delta - (i+j)$, $i+j+k \in 2\Z$}, \\
    0& \text{ otherwise.}
  \end{cases}
  \]
\end{theorem}

At this point we summon a $G$-equivariant fusion category using the algebra in a category construction recalled in \seref{s:algorb}. Suppose that $\delta = 4m$. Put $A= V_0 \oplus V_{\delta}$. Then the $\DD$-algebra $A$ is rigid and $\theta_A = \id$. There is a natural action of the group $\Z_2$ acting by automorphisms $\pi$ on $A$. In particular, the automorphism $\pi_a$ is given by
  \begin{align*}
    \pi &\colon A \to A \text{ where }\\
    &\pi_{a}|_{V_0} = \id \\
    &\pi_{a}|_{V_{\delta}} = -\id
  \end{align*}

  Let $\CC \equiv \rep A$ be the category of modules over $A$. \thref{t:algcat} tells us that this category is a $\Z_2$-equivariant fusion category. \thref{t:modmod} tells us that this category is modular and that the subcategory $\CC_{\e} \subset \CC$ is a modular fusion category. See \firef{f:roadmap1}.

We discussed what the objects in $\CC \equiv \rep A$ are in \seref{s:algorb}. For brevity we refer to representatives of the simple modules by $X_{1}, \ldots , X_{2m-1}$, and $X^{+}_{2m}$ and $X^{-}_{2m}$. The $X_{i}$ for $i \in \{0,1, \ldots, 2m-1\}$ in the subring generated by $X_{1}$ are given by $X_i = V_i \oplus V_{\delta - i} = A \ttt V_i$. The modules $X^{\pm}_{2m}$ are isomorphic as objects of $\DD$ to $V_{2m}$. The object $X_{2m}^{+} \oplus X_{2m}^{-}$ is also in the subring generated by $X_{1}$ and is given by $A \ttt V_{2m}$. See \cite{KO} for more specific details.

Denote the non-trivial element of the group $\Z_{2}$ by $a$. The action of $\Z_2$ on $A$ gives rise to the $\Z_2$-grading and action of $\Z_2$ on $\CC$. The group grading and action on the category is as follows. $X_{i}$ is in $\CC_{\e}$ if and only if $i$ is even, otherwise $X_{i} \in \CC_{a}$. The non-trivial element $a \in \Z_2$ acts trivially on the isomorphism classes of $X_i$ for $i = 1, \dots, 2m-1$ and it interchanges $X^{+}_{2m}$ and $X^{-}_{2m}$.

\cite[Section 7]{KO} gives us the tensor decomposition for the simple objects in the category.

\begin{theorem}\label{t:tensdecomp} Let $\CC$ be the category of representations of the algebra $A$ as described above for $\delta=4m$. Suppose $8|\delta$. Then the decomposition of the tensor product in $\CC$ is given by

  \begin{align*}
  X_{0} \ttt& X_{i} = X_{i}, \\
  X_{1} \ttt& X_{i} = X_{i-1} \oplus X_{i+1}, \mspace{9mu}%
  i=1, \ldots, 2m-2, \\
  X_{1} \ttt& X_{2m-1} = X_{2m-2} \oplus X_{2m}^{+}%
  + X_{2m}^{-}, \\
  X_{1} \ttt& X_{2m}^{\pm} = X_{2m-1}, \\
  X_{2m}^{\pm} \ttt& X_{2m}^{\pm} = X_{0} \oplus X_{4}%
  \oplus \cdots \oplus X_{2m-4} \oplus X_{2m}^{\pm}, \\
  X_{2m}^{\pm} \ttt& X_{2m}^{\mp} = X_{2} \oplus X_{6}%
  \oplus \cdots \oplus X_{2m-2}.   
\end{align*}

The remaining tensor products can be derived from this table. 
\end{theorem}

The results for $\delta \equiv 4 \mod 8$ are similar and left omitted. In this entire section we shall consider only the case where $\delta=4m$ and $8|\delta$.

 \subsection{Example} Suppose we wish to find $X_2 \ttt X^{\pm}_{2m}$. Since $X_1 \ttt X_1 = X_0 \oplus X_2$ consider 
  \begin{align*}
    X_1 \ttt X_1 \ttt X^{\pm}_{2m} &= \left( X_0 \oplus X_2 \right) \ttt X^{\pm}_{2m}\\
    &= \left( X_0 \ttt X^{\pm}_{2m} \right) \oplus \left( X_2 \ttt X^{\pm}_{2m} \right)\\
    &= X^{\pm}_{2m} \oplus \left( X_2 \ttt X^{\pm}_{2m} \right). \\
    \intertext{And } 
    X_1 \ttt X_1 \ttt X^{\pm}_{2m} &= X_1 \ttt  X_{2m-1}  \\
    &= X_{2m-2} \oplus X^{+}_{2m} \oplus X^{-}_{2m}. \\
    \intertext{So that }
    X_2 \ttt X^{\pm}_{2m} &= X_{2m-2} \oplus X^{\mp}_{2m}.
  \end{align*}

The notion of a $G$-equivariant fusion category was not used in \cite{KO} and the extended Verlinde algebra had not yet been formalized. Regardless, although in a different form, the coefficients for fusion in $\Ve_{\e} \subset \Ve$ are given above by \thref{t:tensdecomp}. We now put those results into the language of the extended Verlinde algebra.

Consider the extended Verlinde algebra $\Ve(\CC)$. Fix a basis and notation for this algebra as described in \seref{s:basis}. Recall, in particular, that $\lambda_{i}$ is the identity map of $X_{i}$. A basis for the sub-algebra $\Ve_{\e} \subset \Ve$ is given by $\lambda_{i}$ where $i \in I$ and $I = I(\CC) = \{0,1, \ldots ,2m-1, {}^{+}_{2m}, {}^{-}_{2m} \}$. Note that $I = I_{\e} \sqcup I_{a}$ where $I_{\e}$ contains the even integers and $I_{a}$ the odd. To refer to the element $\lambda_{2m}^{+} + \lambda_{2m}^{-}$ we sometimes will use the notation $\lambda_{2m} = (\lambda_{2m}^{+} + \lambda_{2m}^{-})$. We can then use the index set $I^{\circ} = \{0,1, \ldots, 2m-1, 2m \}$ to refer to the ring generated by $\lambda_{1}$. Then the fusion rules for elements in $\Ve_{\e}$ are given by the above theorem where we replace $X_{i}$ by $\lambda_{i}$.

The results from \cite{KO} however do not give a complete description of the fusion rules for this algebra. In particular the fusion with elements of the sort $\basis[a]{i}$ are yet to be considered. Later,  \seref{s:Z2}, we will consider such elements when we show that the fusion rules can be diagonalized in the extended Verlinde algebra arising from any modular $\Z_{2}$-equivariant fusion category.

Right now we want to check that the formula given in \coref{c:v11mod} correctly predicts the fusion rules. To use this formula we need to have the operator $s^{\CC}$ at our disposal. The operator $s^{\CC}$ can be given in terms of the operator $s^{\DD}$ except for the elements $\lambda^{+}_{2m}$ and $\lambda^{-}_{2m}$. Let us christen these elements, the exceptional elements; and give them the notation $\lambda^{+}$ and $\lambda^{-}$ respectively. Recall that it is on exactly the exceptional elements which the group acts non-trivially. To find $s^{\CC}$ on the exceptional elements we will have to refer to the universal $R$-matrix of $\DD$. We proceed in two parts. We first discuss the example on the non-exceptional elements. We give the relation between the $s$-matrix of the two algebras, the relation between the fusion coefficients, and we show that the corollary does indeed produce the fusion rules. We dispatch the second part by using a symmetry argument and the result from the first: the tensor product decomposes in a similar way for similar products, the $s$-matrix has the same value on similar points. Afterwards, we explicitly find the value of the $s$-matrix at some exceptional points so that one can use the standard Verlinde formula in the subalgebra $\Ve_{\e,\e} \subset \Ve_{\e}$.

Recall that the $s^{\DD}$-matrix for $\U_q(\mathfrak{sl}_2)$ is given by the following

\begin{theorem}\label{t:smatver}
  \[
  (s \chi_{i}, \chi_{j}) = \sqrt{\frac{2}{\varkappa}} \sin \left(\frac{(i+1)(j+1)\pi}{\varkappa}\right).
  \]
\end{theorem}

\subsection{Example}\thref{t:srelat} can be tricky to work with in practice. So we make a preemptive strike by calculating two examples explicitly. We shall find $(s \lambda_2 , \lambda_2)$ and $(s \lambda_3, \basis[a]{2})$ in terms of $s^{\DD}$ using \coref{c:srelat}.

Let $x=\chi_{2} \in \V(\DD)$ and $y=\lambda_{2}$. Then 
  \begin{align*}
   (s \widehat F' x,y)_{\CC} &= (s \lambda_{2}, \lambda_{2})_{\CC} + (s \basis[a]{2}, \lambda_{2})_{\CC}\\
    &= (s \lambda_{2}, \lambda_{2})_{\CC}
  \end{align*}
  since by definition $(sP,Q) = 0$ if $P \in \Ve_{a,b}$ but $Q \notin \Ve_{b,a}$. And  \begin{align*}
   (s x, \widehat G'y)_{\CC / G} &= (s \chi_{2},\chi_{2})_{\CC / G} + (s \chi_{2}, \chi_{2})_{\CC / G}\\
    &= 2 (s \chi_{2}, \chi_{2})_{\CC / G}.
  \end{align*}
  So that, $(s \lambda_{2}, \lambda_{2})_{\CC} = 2 (s \chi_{2}, \chi_{2})_{\CC / G}$.

Now let $x=\chi_{2} \in \V(\DD)$ and $y=\lambda_{3}$. Then $(s \widehat F' x,y)_{\CC}$
  \begin{align*}
  &= (s \lambda_{2}, \lambda_{3})_{\CC} + (s \basis[a]{2}, \lambda_{3})_{\CC}\\
  &= (s \basis[a]{2}, \lambda_{3})_{\CC}.
  \end{align*}
And $(s x, \widehat G'y)_{\CC / G}$
  \begin{align*}
    &= (s \chi_{2},\chi_{3})_{\CC / G} + (s \chi_{2}, \chi_{3})_{\CC / G}\\
    &= 2 (s \chi_{2}, \chi_{3})_{\CC / G}.
  \end{align*}
  So that, $(s \lambda_{3}, \basis[a]{2})_{\CC} = 2 (s \chi_{3}, \chi_{2})_{\CC / G}$.

By the equivalence $\DD \backsimeq \CC / G$ we have 
  \begin{align*}
    (s \lambda_{2}, \lambda_{2})_{\CC} &= 2 (s \chi_{2}, \chi_{2})_{\DD} \\
    \text{and } (s \lambda_{3}, \basis[a]{2})_{\CC} &= 2 (s \chi_{3}, \chi_{2})_{\DD}.
  \end{align*}

  \vspace{10pt}
Rather than prove by hand that the $s^{\CC}$-matrix and the corollary give the fusion rules we shall instead establish the result by using \thref{t:quantcg} and the functor $F$ introduced earlier. This turns out to be a surprisingly good way to proceed. We observe at the level of the Verlinde formula what happens to the symmetry in $\DD$ when we pass from the fusion category $\DD$ to the equivariant fusion category $\CC$.

In the sequel we will need a lemma. In order not to lose focus we avoid giving it in a more general setting. Consider it a technical lemma.

\begin{lemma}\label{l:tech}
  Let $F \colon \DD \to \CC$ and $G \colon \CC \to \DD$ be the adjoint functors given in \seref{s:algorb}. If $i \in I^{\circ}$ then
  \[ F(G(\lambda_{i})) = 2 \lambda_{i}.\]
\end{lemma}

\begin{proof}
  Recall that $\lambda_{i} \colon X_{i} \to X_{i}$ is an identity map. Since $i \in I^{\circ}$ we have that 
  \[A \ttt X_i = (V_0 \oplus V_{\delta}) \ttt X = X \oplus X \]
  and the lemma.
\end{proof}

To show that the corollary gives the fusion rules we first re-cast \thref{t:tensdecomp}.

\begin{theorem}\label{t:coefrelat}
Let $\CC \equiv \rep A$ as described above for $\delta=4m$. Suppose $i,j \in I^{\circ}$. Write \( \lambda_{i} \ttt \lambda_{j} = \sum_{k} L_{ij}^{k} \lambda_{k} \) where $k \in I^{\circ}$. Then 
  \[
  L_{ij}^{k} =\begin{cases}
    N_{ij}^{k} + N_{ij}^{\delta - k}& \mspace{18mu} \text{for $k \neq 2m$,} \\
    N_{ij}^{k}& \mspace{18mu} \text{for $k = 2m$,}
  \end{cases}
  \]
where $N_{ij}^{k}$ are the fusion coefficients in $\V(\DD)$.
\end{theorem}

\begin{proof}
  The element $\lambda_{i} \in \Ve$ is the identity map $\lambda \colon X_{i} \to X_{i}$. By \leref{l:tech} we can write $\lambda_{i} \ttt \lambda_{j}$ as
  \[
  = \frac{1}{2} F(G(\lambda_{i} \ttt \lambda_{j})) = \frac{1}{2} F( \chi_{i} \ttt \chi_{j} \ttt (\chi_{0} + \chi_{\delta})) 
  \]
  where $\chi_{i} \in \V(\DD)$ is the identity map $\chi \colon V_i \to V_i$. Now $\chi_{i} \ttt \chi_{j}$ can be written as $\sum_{k} N_{ij}^{k} \mspace{3mu} \chi_{k}$, where $k \in I(\DD)$. So $\lambda_{i} \ttt \lambda_{j}$ can be written as follows
  \[
  \begin{split}
  &= \frac{1}{2} F \sum_{k \neq 2m} N_{ij}^{k} \mspace{3mu} \chi_{k} \ttt \left( \chi_{0} + \chi_{\delta} \right) \mspace{9mu}+\mspace{9mu} N_{ij}^{2m} \mspace{3mu} (\chi_{2m} \ttt (\chi_{0} + \chi_{\delta})) \\
  &= \frac{1}{2} F \sum_{k \neq 2m} N_{ij}^{k} \left( (\chi_{k} \ttt \chi_{0}) + (\chi_{k} \ttt \chi_{\delta}) \right) \mspace{9mu}+\mspace{9mu} N_{ij}^{2m} \mspace{3mu} (\chi_{2m} \ttt (\chi_{0} + \chi_{\delta})) \\
  &= \frac{1}{2} F \sum_{k \neq 2m} N_{ij}^{k} \mspace{3mu} (\chi_{k} + \chi_{\delta - k}) \mspace{9mu}+\mspace{9mu} N_{ij}^{2m} \mspace{3mu} (\chi_{2m} \ttt (\chi_{0} + \chi_{\delta}))  \\
  &= \frac{1}{2} F \sum_{k = 0}^{2m-1} \mspace{-6mu} N_{ij}^{k} \mspace{3mu} (\chi_{k} + \chi_{\delta - k}) \mspace{6mu} + \mspace{-9mu} \sum_{k = 2m+1}^{\delta} \mspace{-15mu} N_{ij}^{k} \mspace{3mu} (\chi_{k} + \chi_{\delta - k}) \\
  & \mspace{36mu} + N_{ij}^{2m} \mspace{3mu} (\chi_{2m} \ttt (\chi_{0} + \chi_{\delta}))  \\
  &= \frac{1}{2} F \sum_{k = 0}^{2m-1} \mspace{-6mu} N_{ij}^{k} \mspace{3mu} (\chi_{k} + \chi_{\delta - k}) \mspace{6mu} + \mspace{-6mu} \sum_{k = 0}^{2m-1} \mspace{-6mu} N_{ij}^{\delta - k} \mspace{3mu} (\chi_{k} + \chi_{\delta - k}) \\
  & \mspace{36mu} + N_{ij}^{2m} \mspace{3mu} (\chi_{2m} \ttt (\chi_{0} + \chi_{\delta})) \\
  &= \frac{1}{2} F \sum_{k = 0}^{2m-1} \mspace{-6mu} (N_{ij}^{k} + N_{ij}^{\delta - k}) \mspace{3mu} (\chi_{k} + \chi_{\delta - k}) \mspace{6mu} + N_{ij}^{2m} \mspace{3mu} (\chi_{2m} \ttt (\chi_{0} + \chi_{\delta})) \\
  &= \sum_{k = 0}^{2m-1} \mspace{-6mu} (N_{ij}^{k} + N_{ij}^{\delta - k}) \mspace{3mu} \lambda_{k} \mspace{6mu} + N_{ij}^{2m} \mspace{3mu} \lambda_{2m} \\
  \end{split}
  \]
\end{proof}

\subsection{Example} Suppose $\delta=8$. In $\V(\DD)$ we have $\chi_{2} \ttt \chi_{3} = \chi_{1} + \chi_{3} + \chi_{5}$. In $\Ve(\CC)$ we have $\lambda_{2} \ttt \lambda_{3} = \frac{1}{2} F(G(\lambda_{2} \ttt \lambda_{3}))$
  \[
  \begin{split}
  &= \frac{1}{2} F\left( \chi_{2} \ttt \chi_{3} \ttt ( \chi_{0} + \chi_{8} ) \right) \\
  &= \frac{1}{2} F\left( \chi_{1} + \chi_{3} + \chi_{5} \ttt (\chi_{0} + \chi_{8}) \right) \\
  &= \frac{1}{2} F\left( (\chi_{1} + \chi_{3} + \chi_{5}) + (\chi_{7} + \chi_{5} + \chi_{3}) \right) \\
  &= \lambda_{1} + 2 \lambda_{3}
  \end{split}
  \]

To show that \coref{c:v11mod} correctly predicts the fusion coefficients we need to prove the following theorem which relates it to the Verlinde formula.

\begin{theorem}\label{t:twosums}
Suppose $i,j,k \in I^{\circ}$. Suppose further that $i$ is even. Then
  \[
  \begin{split}
  &\sum_{p} \frac{(s \lambda_i, \lambda_p)(s \lambda_j, \basis{p})(s \lambda_k^*, \basis{p})}{(s \lambda_0, \lambda_p)}\\
  &\mspace{27mu}=
  \left\{
  \begin{aligned}
  &\sum_{p'} \, \frac{(s \chi_{i} , \chi_{p'})(s \chi_{j} , \chi_{p'}) \left((s \chi_{k}^{*} , \chi_{p'})+(s \chi_{\delta - k}^{*} , \chi_{p'}) \right)}{(s \chi_0 , \chi_{p'})} \mspace{9mu} &\text{for } k \neq 2m \\
  &\sum_{p'} \, \frac{(s \chi_{i} , \chi_{p'})(s \chi_{j} , \chi_{p'}) (s \chi_{k}^{*} , \chi_{p'})}{(s \chi_0 , \chi_{p'})}  &\text{for } k = 2m \\
  \end{aligned}
  \right.
  \end{split}
  \]
where $p$ is indexed by $I_{\e} \subset I(\CC)$ and $p'$ is indexed by $I' = I(\DD)$.
\end{theorem}

To prove the theorem we first establish several interesting facts. The symmetry of the tensor product in $\DD = \rep \U_{q}(\mathfrak{sl}_{2})$ which is in some sense a structure {\em on} the Verlinde algebra $\V(\DD)$ becomes part of the structure {\em in} the extended Verlinde algebra $\Ve(\CC) = \Ve(\rep A)$.

\begin{lemma}\label{l:oddfold} Suppose that $p'$ is odd. Then
  \[
  (s \chi_{k}^{*} , \chi_{p'})+(s \chi_{\delta - k}^{*} , \chi_{p'}) = 0.
  \]
\end{lemma}

\begin{corollary}\label{c:oddfold} Suppose that $p'$ is odd and that $k=2m$. Then
  \[
  (s \chi_{k}^{*}, \chi_{p'}) = 0.
  \]
\end{corollary}

\begin{lemma}\label{l:evenfold} Suppose that $p'$ is even. Then
  \[
  (s \chi_{k}^{*} , \chi_{p'}) = (s \chi_{\delta - k}^{*} , \chi_{p'}).
  \]
\end{lemma}

These facts follow essentially because the sine function is odd with respect to reflection about $2\pi$, and even with respect to reflection about $\pi/2$. Note that these results will hold in the more general case when passing from a modular fusion category to a modular $G$-equivariant fusion category via the algebra in a category construction.

\begin{proof}[Proof of \leref{l:oddfold}]Note that $(\delta - k + 1)(2w)$
  \[
  \begin{split}
  &= 2w\delta - 2wk + 2w + (2w -2w) \\
  &= 2w\delta + 4w - (2wk+2w) \\
  &= (2w)(\delta +2) - (k+1)(2w).
  \end{split}
  \]
  So that
  \[
  \begin{split}
  (s \chi_{\delta - k}^{*}, \chi_{2w-1}) &= \sqrt{2/{(\delta+2)}} \sin \left( \frac{(\delta - k + 1)(2w)\pi}{(\delta + 2)} \right) \\
  &= \sqrt{2/{(\delta+2)}} \sin \left( \frac{-(k+1)(2w)\pi}{(\delta + 2)} \right) \\
  &= {-1} \cdot \sqrt{2/{(\delta+2)}} \sin \left( \frac{(k + 1)(2w)\pi}{(\delta + 2)} \right) \\
  &= {-1} \cdot (s \chi_{k}^{*}, \chi_{2w-1})
  \end{split}
  \]
  and the lemma.
\end{proof}

\begin{proof}[Proof of \coref{c:oddfold}]The lemma implies that
  \[
  (s \chi_{2m}, \chi_{p'}) + (s \chi_{\delta-2m}, \chi_{p'})
  \]
is zero. Note that $\delta - 2m = 4m - 2m = 2m$ and the result follows.
\end{proof}

\begin{proof}[Proof of \leref{l:evenfold}] Note that $(\delta - k + 1)(2w+1)$
  \[
  \begin{split}
  &= 2w\delta + \delta - 2wk - k + 2w + 1 + ((2w+1) - (2w+1)) \\
  &= 2w\delta + 4w + \delta + 2 - 2wk - k - (2w+1) \\
  &= 2w(\delta +2) + (\delta +2) - (k(2w+1) + (2w+1)).
  \end{split}
  \]
The proof follows the previous lemma with the nuance that here we reflect across $\pi/2$.
\end{proof}

\begin{proof}[Proof of \thref{t:twosums}]
Consider the right-hand side of the equation when $k \neq 2m$. Here we are summing over $p' \in I(\DD)$. By \leref{l:oddfold} we have that the summand is zero when $p'$ is odd. When $j$ is odd, that is $\lambda_{j} \in \Ve_{\e,a}$, we have that the summand corresponding to $p' = 2m$ is zero by \coref{c:oddfold}. So, in this case, the sum over $I$ is equivalent to the sum over $\{0,2,\ldots,2m-2,2m+2,\ldots,\delta\}$. Now apply \leref{l:evenfold} and \leref{l:oddfold}. The right-hand side simplifies thus
  \[
  4 \cdot \sum_{p'} \, \frac{(s \chi_{i} , \chi_{p'})(s \chi_{j} , \chi_{p'}) (s \chi_{k}^{*} , \chi_{p'})}{(s \chi_0 , \chi_{p'})},
  \]
where $p'$ runs over the set $\{0,2,\ldots,2m-2\}$. In the case that $\lambda_{j} \in \Ve_{\e,\e}$ \coref{c:oddfold} does not apply, so we get a $2m$ summand. Since $j$ is even we apply \leref{l:evenfold} twice to get the same indexing set as in the odd case. The right-hand side is written thus
  \[
  4 \cdot \sum_{p'} \, \frac{(s \chi_{i} , \chi_{p'})(s \chi_{j} , \chi_{p'}) (s \chi_{k}^{*} , \chi_{p'})}{(s \chi_0 , \chi_{p'})} \mspace{9mu}+\mspace{9mu}  2 \cdot \frac{(s \chi_{i} , \chi_{2m})(s \chi_{j} , \chi_{2m}) (s \chi_{k}^{*} , \chi_{2m})}{(s \chi_0 , \chi_{2m})} ,
  \]
where $p'$ runs over the same indexing set. Now consider the left-hand side of the equation. $p$ is indexed by $I_{\e}$. Recall that $I_{\e} = \{0,2,\ldots,2m-2,{}_{2m}^{+},{}_{2m}^{-}\}$. Use the relation between $s^{\DD}$ and $s^{\CC}$ given by \coref{c:srelat} to get equality between the two sides. We note that the case $k=2m$ runs similar and conclude the proof.

\end{proof}

\subsection{Exceptional Elements} Now we turn our attention to the exceptional elements in the algebra. We want to show that the fusion rules for the exceptional elements are also correctly predicted by \coref{c:v11mod}. The exceptional elements $\lambda^{+}$ and $\lambda^{-}$ are in the subalgebra $\Ve_{\e,\e}(\CC) \subset \Ve(\CC)$. Recall \firef{f:roadmap1}. Suppose that $\lambda_{j} \in \Ve_{\e,a}$. Then we know {\it a priori} that $\lambda^{\pm} \ttt \lambda_{j}$ belongs to $\Ve_{\e,a}$. That is, the number $L_{\pm,j}^{k}$ is nonzero only when $k \in I_{a}$. Note that since the dimension of $\lambda^{+}$ and $\lambda^{-}$ are equal we have that: for all $p \in {}^{a}I_{\e}$ that $(s \lambda^{+},\lambda_{p}) = (s \lambda^{-},\lambda_{p})$. We have already shown that \coref{c:v11mod} holds in the case that $i=2m$ and $j,k \in I_{a}$. Recall that $L_{2m,j}^{k}= L_{+,j}^{k} + L_{-,j}^{k}$. Then since $\lambda^{+} \ttt \lambda_{j} = \lambda^{-} \ttt \lambda_{j}$, stated differently, that $L_{+,j}$ is identical to $L_{-,j}$,  we have that \coref{c:v11mod} correctly predicts the fusion rules for this example. Indeed,

  \[
  \begin{split}
    2 L_{\pm,j}^{k} &= L_{2m,j}^k = \sum_{p} \frac{(s \lambda_{2m}, \lambda_p)(s \lambda_j, \basis{p})(s \lambda_k^*, \basis{p})}{(s \lambda_0, \lambda_p)} \\
    &= \sum_{p} \frac{\left( (s \lambda^{+}, \lambda_p) + (s \lambda^{-}, \lambda_p) \right) (s \lambda_j, \basis{p})(s \lambda_k^*, \basis{p})}{(s \lambda_0, \lambda_p)}\\
    &= 2 \sum_{p} \frac{(s \lambda^{\pm}, \lambda_p)(s \lambda_j, \basis{p})(s \lambda_k^*, \basis{p})}{(s \lambda_0, \lambda_p)}.\\
    \end{split}
    \]

We have not considered yet the fusion of two objects in $\Ve_{\e,\e}$ when one of them is an exceptional element. We do this now. Again, the subalgebra $\Ve_{\e,\e} \subset \Ve_{\e}$ is a standard Verlinde algebra. We have the usual Verlinde formula at our disposal to calculate the fusion rules. The only requisite to using the Verlinde formula is having the values of the $s$-matrix in the Verlinde algebra. \thref{t:srelat}, remarkably enough, can be used to find the values of the $s$-matrix in $\Ve_{\e}(\CC)$ except for the entries
  
  \[
  \begin{split}
  &(s \lambda_{2m}^{\pm},\lambda_{2m}^{\pm}) \\
  \text{and }  &(s \lambda_{2m}^{\pm},\lambda_{2m}^{\mp}).
  \end{split}
  \]

There are several ways to find these values. Let us restrict the discussion to the case where $\delta = 4m$ and $8|\delta$. Then the product in $\V(\DD)$
  \begin{align*}
    \chi_{2m} \ttt \chi_{2m} &= \chi_{0} + \chi_{2} + \ldots + \chi_{2m} + \ldots + \chi_{4m-2} + \chi_{4m} \\
    \intertext{splits in $\Ve(\CC)$ as follows}
    \lambda_{2m}^{\pm} \ttt \lambda_{2m}^{\pm} &= \lambda_{0} + \lambda_{4} + \ldots + \lambda_{2m-4} + \lambda_{2m}^{\pm} \mspace{18mu} \text{and}\\
    \lambda_{2m}^{\pm} \ttt \lambda_{2m}^{\mp} &= \lambda_{2} + \lambda_{6} + \ldots + \lambda_{2m-2}.
  \end{align*}

This gives us the immediate relation
\begin{lemma}\label{l:excsum}
  \[
  (s^{\DD} \chi_{2m},\chi_{2m}) = (s^{\CC} \lambda_{2m}^{\pm},\lambda_{2m}^{\pm}) + (s^{\CC} \lambda_{2m}^{\pm},\lambda_{2m}^{\mp})
  \]
\end{lemma}

To find the values $(s \lambda_{2m}^{\pm},\lambda_{2m}^{\pm})$ and $(s \lambda_{2m}^{\pm},\lambda_{2m}^{\mp})$ we refer to the universal $\widecheck{R}$-matrix of $\DD \equiv \rep \U_q(\mathfrak{sl}_2)$ . In $\V(\DD)$ the relevant formula is given by
  \begin{gather*}
  \raisebox{-0.5\height}{\begin{overpic}
    {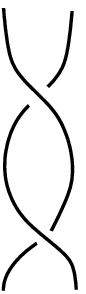}
    \put(-12,5){$\chi_{i}$}
    \put(30,5){$\chi_{j}$}
  \end{overpic}}
  \mspace{21mu}=\mspace{21mu}\sum_{k} \mspace{3mu} \widecheck{R}^{2}|_{k} \mspace{3mu} N_{ij}^{k}
\raisebox{-0.5\height}{\begin{overpic}
    {arrow}
    \put(10,5){$\chi_{k}$}
  \end{overpic}} \\
  \intertext{where $\widecheck{R}^{2}$ is given by}
  \widecheck{R}^{2}|_{V_{k} \subset V_{i} \ttt V_{j}} = q^{-\frac{1}{2}(i(i+2) + j(j+2) - k(k+2))} = \theta_{i}^{-1}\theta_{j}^{-1}\theta_{k}.
\end{gather*}

To get the value $(\tilde s \chi_{i},\chi_{j})$ take the trace of both sides. That is, close the strands.

\subsection{Example} Suppose that $\delta = 8$. Let us find $(s \lambda_{4}^{\pm}, \lambda_{4}^{\pm})$ and $(s \lambda_{4}^{\pm}, \lambda_{4}^{\mp})$. Using the functors $F$ and $G$, and the relation $dim_{\CC}(X) = \frac{dim_{\DD}(X)}{dim_{\DD}A}$ from \thref{t:dimrel} yields the following equation

  \[
  \raisebox{-0.5\height}{\begin{overpic}
  {smatrix10}
  \put(3,58){$\lambda^{\pm}$}
  \put(49,11){$\lambda^{\pm}$}
  \end{overpic}}
  \mspace{9mu}=\mspace{9mu}\frac{q^{-4}}{2}\mspace{3mu}
  \raisebox{-0.5\height}{\begin{overpic}
  {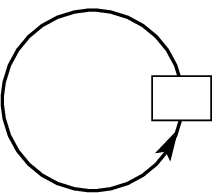}
  \put(77,42){$\chi_{0}$}
\end{overpic}}
\mspace{9mu}+\mspace{9mu}\frac{q^{8}}{2}\mspace{3mu}
 \raisebox{-0.5\height}{\begin{overpic}
{smatrix4}
\put(77,42){$\chi_{4}$}
\end{overpic}}
\mspace{9mu}+\mspace{9mu}\frac{q^{-4}}{2}\mspace{3mu}
 \raisebox{-0.5\height}{\begin{overpic}
{smatrix4}
\put(77,42){$\chi_{8}$}
\end{overpic}}
 \]

Thus $(\tilde s \lambda_{4}^{\pm}, \lambda_{4}^{\pm}) = \frac{1}{2}(2 q^{-4} [1] + q^8 [5])$. After normalization and \leref{l:excsum} we get 
\begin{align*}
 (s \lambda_{4}^{\pm}, \lambda_{4}^{\pm}) &= \sqrt{\frac{2}{10}} \left(1- 2 \sin\frac{3 \pi}{10} \right), \\
 (s \lambda_{4}^{\pm}, \lambda_{4}^{\mp}) &= \sqrt{\frac{2}{10}} \left(2 \sin\frac{3 \pi}{10}\right).
 \end{align*}

 \begin{theorem}\label{t:excval} Let $\delta=4m$ and $\varkappa = \delta+2$. Suppose that $8|\delta$. Then
  \[
  (s \lambda^{\pm}_{2m}, \lambda^{\pm}_{2m}) = \frac{1}{2} \left(\sqrt{\frac{2}{\varkappa}} + (-1)^{m/2}\right).
  \]
\end{theorem}

\vspace{5pt}
Before giving the proof of the theorem we recall [See \thref{t:smatver}], for the sake of comparison and curiosity, that
\[(s \chi_{2m},\chi_{2m})^{\DD}  = \sqrt{\frac{2}{\varkappa}}. \]
\vspace{5pt}

\begin{proof}
 Apply the formula $\tilde s_{ij} = \theta^{-1}_{i}\theta^{-1}_{j} \sum_{k} N_{i^{*}j}^{k}\theta_{k} d_{k}$ to the decomposition of $\lambda_{2m}^{\pm} \ttt \lambda_{2m}^{\pm}$ to get
 \begin{align*}
   (\tilde s \lambda^{\pm}, \lambda^{\pm}) &= \theta_{2m}^{-2} \sum_{p=0}^{m} \theta_{4p} dim_{\CC}(\lambda_{4p}) = \theta_{2m}^{-2} \frac{1}{dim_{\DD}A} \sum_{p=0}^{m} \theta_{4p} [4p+1] \\
   &= \frac{\theta_{2m}^{-2}}{q^{1} - q^{-1}} \frac{1}{2} \sum_{p=0}^{m} q^{{1/2}(4p)(4p+2)} (q^{4p+1} - q^{-1(4p+1)}) \\ 
   &= \frac{\theta_{2m}^{-2}q^{-1}}{q^{1} - q^{-1}} \frac{1}{2} \sum_{p=0}^{m} q^{2(2p+1)^{2}} - q^{2(2p)^{2}} \\
   &= \frac{\theta_{2m}^{-2}q^{-1}}{q^{1} - q^{-1}} \frac{1}{2} \sum_{p=0}^{2m+1} {(-1)}^{p+1} q^{2p^{2}}.
\end{align*}

Recall that $q = e^{\pi i / \varkappa}$. Use the symmetries $-q^{2p^{2}} = q^{2(p+\varkappa/2)^{2}}$ and $q^{2p^{2}} = q^{2(p+\varkappa)^{2}}$ to further simplify the sum:
  \[
  (\tilde s \lambda^{\pm}, \lambda^{\pm}) = \frac{\theta_{2m}^{-2}q^{-1}(-1)}{q^{1} - q^{-1}} \frac{1}{2} \left(1 + \frac{1}{2}\sum_{p=1}^{\varkappa} q^{8p^{2}}\right).
  \]

  Explicit calculation shows that $-\theta_{2m}^{-2}q^{-1} = -q^{-(2m+1)^{2}} = - (-i)^{2m+1} = (-1)^m(i)$ and $q^{1}-q^{-1} = 2i \sin (\pi / \varkappa)$ so that the coefficient 
  \[
  \frac{\theta_{2m}^{-2}q^{-1}(-1)}{q^{1} - q^{-1}} \frac{1}{2} = \frac{(-1)^m}{4 \sin(\pi /\varkappa)}.
  \]

  The sum $\sum_{p=1}^{\varkappa} q^{8p^{2}} = \sum_{p=1}^{\varkappa} e^{8 \pi i p^{2} / \varkappa}$ is a quadratic Gauss sum. Following the convention in \cite{Ap} let
  \[
  S(a,b) = \sum_{p=1}^{b} e^{\pi i a p^{2}/b}.
  \]
  The sums $S(a,b)$ enjoy a ``reciprocity law'' (see, e.g., \cite[Section 9.10]{Ap}): if $ab$ is even then
  \[
  S(a,b) = \sqrt{\frac{b}{a}}\left(\frac{1+i}{\sqrt{2}}\right)\overline{S(b,a)}.
  \]
It is straight-forward enough to find $S(\varkappa,8) = 2\sqrt{2}(1+i)i^{m}$ which gives $S(8,\varkappa) = \sqrt{2\varkappa}(-i)^{m}$ and 
  \[
  (\tilde s \lambda^{\pm}, \lambda^{\pm}) = \frac{(-1)^{m}}{4 \sin(\pi /\varkappa)} \cdot \left(1+(-i)^{m}\sqrt{\frac{\varkappa}{2}}\mspace{6mu}\right).
  \]
  To conclude the proof note that $m$ is even by the hypotheses and normalize by multiplying by $D^{-1}_{\CC} = \left|G\right| \sqrt{2/\varkappa}\sin(\pi/\varkappa)$.
\end{proof}

\begin{corollary}
  \[
  (s \lambda^{\pm},\lambda^{\mp}) = \frac{1}{2}\left(\sqrt{\frac{2}{\varkappa}} - (-1)^{m/2} \right).
  \]
\end{corollary}

\begin{proof}
  Immediate from \thref{t:excval} and the relation given by \leref{l:excsum}.
\end{proof}

In the next section we extend our results to fusion when our category is a general modular $\mathbb{Z}_2$-equivariant fusion category. At the end of the next section we will consider the example $D_{2m+2}$ again, where the case is that both $\lambda_{i}$ and $\lambda_{j}$  are in $\Ve_{\e,a}$.

\section{Fusion in $\Ve$ where $G=\mathbb{Z}_2$}\label{s:Z2}

In \seref{s:uqsl} we showed that the extended Verlinde formula correctly predicts the fusion rules of the extended Verlinde algebra arising from the algebra of type $D_{2m+2}$. The category came equipped with a $\mathbb{Z}_2$ grading and action. This was the first time that we considered the fusion of two elements that were not both in the space $\Ve_{\e,\e}$. In this section we show that we can diagonalize the fusion rules in the general case where $G = \mathbb{Z}_2$. 

For this section we suppose that $\CC$ is a modular $\mathbb{Z}_2$-equivariant category. We denote the element in $\mathbb{Z}_2$ which generates the group by $a$. Again, for $\ph$ and $\psi$ in $\Ve_{\e,a}$ we are interested in knowing how $\ph \ttt \psi$ decomposes. From the definition of the tensor product we know a priori that this tensor will decompose as a sum in $\Ve_{\e,\e}$. This is a significant change from the setting in \seref{s:v11} where the operator $L_\ph$ [which was left multiplication by $\ph \in \Ve_{\e,\e}$] took the space $\Ve_{\e,g}$ to itself.

In \seref{s:v11} we obtained our results by using the $s$-matrix to interchange the tensor product with the convolution product. For $\ph \in \Ve_{\e,\e}$ we had that $\tilde s \ph$ was a linear combination of identity maps. This made convolving an easy process. We no longer have this nicety. We need another approach. Roughly speaking, we will pick a generic basis and then find a change of basis which diagonalizes the convolution product. Let us make this formal now.

\subsection{Diagonalizing the Convolution Product}\label{s:diag}

Consider the space $\Ve_{*,\e}$. This space is $\Ve_{\e,\e} \oplus \Ve_{a,\e}$. As described in \seref{s:basis} we choose a generic basis for $\Ve$ with a slight change. First recall that for the space $\Ve_{\e,*}$ we picked the identity morphisms as the basis elements; these identity morphisms were denoted by $\lambda_{i}$. For $i \in {}^aI_{\e}$ pick basis elements $\basis[a]{i}$ for the space $\Ve_{a,\e}$ with the following property

	\[
	\basis[a]{i} \mspace{3mu} * \basis[a]{i} \mspace{3mu} = \mspace{3mu} d_i^{-1} \mspace{3mu} \lambda_{i}.
	\]
	
We can obviously pick such a $\basis[a]{i}$. Indeed, for $i \in {}^aI_{\e}$ and a fixed representative $V_i$ in $\CC$, we can pick an arbitrary element $\ph \colon V_i \to {}^aV_i$ from the space $\Ve_{a,\e}$. Then $\ph * \ph$ is some [non-zero] multiple of $\lambda_i$. Make the appropriate normalization to get $\basis[a]{i}$.

Fix this basis. We are now prepared to find a suitable change of basis so that convolution becomes diagonalized. For $i \in {}^aI_{\e}$ we define

\begin{align*}
  \alpha_{i} &= \tfrac{1}{2}(\basis[a]{i} - \lambda_{i}), \\
  \beta_{i}  &= \tfrac{1}{2}(\basis[a]{i} + \lambda_{i}). \\
\end{align*}

Define the change of basis operator $M \colon \Ve_{*,\e} \to \Ve_{*,\e}$ as follows

\begin{align*}
  M \lambda_{i} &= \alpha_{i}, \mspace{6mu} i \in {}^aI_{\e} \\
  M \basis[a]{i} &= \beta_{i}, \mspace{6mu} i \in {}^aI_{\e} \\
  M \lambda_{i} &= \lambda_{i}, \mspace{6mu} i \in I_{\e} \mspace{6mu} \textrm{and} \mspace{6mu} i \notin {}^aI_{\e}
\end{align*}

Order the basis so that $\{ \ldots, \mspace{3mu} \lambda_{i} \mspace{3mu}, \basis[a]{i} \mspace{3mu}, \mspace{3mu} \ldots \}$ appear pairwise. Then the change of basis matrix \(M\) will then be a diagonal block matrix. The blocks are

\[
\left[
\begin{smallmatrix}	
  \text{-}\frac{1}{2}& \frac{1}{2}\\
  \text{ }\frac{1}{2}& \frac{1}{2}\\
\end{smallmatrix}
\right]
\mspace{18mu} \text{and} \mspace{18mu}
\left[
\begin{smallmatrix}
  1
\end{smallmatrix}
\right].
\]

A simple calculation shows that the convolution product is diagonal in this basis:

\begin{align*}
  \alpha_{i} * \alpha_{i} \mspace{6mu} =& \mspace{6mu} -d_i^{-1} \alpha_{i} \\
  \alpha_{i} * \beta_{i} \mspace{6mu} =& \mspace{6mu} 0 \\
  \intertext{and}
  \beta_{i} * \beta_{i} \mspace{6mu} =& \mspace{6mu} d_i^{-1} \beta_{i} \\
  \beta_{i} * \alpha_{i} \mspace{6mu} =& \mspace{6mu} 0 \\
  \intertext{and for $i \notin {}^aI_{\e}$}
  \lambda_{i} * \lambda_{i} \mspace{6mu} =& \mspace{6mu} d_i^{-1} \lambda_{i}.
\end{align*}

Convolution in $\Ve_{*,\e}$ is now diagonal. This approach will allow us to find a nice equation for the fusion rules. Suppose that $\ph$ and $\psi$ are in $\Ve_{\e,a}$. Then $\tilde s \ph$ and $\tilde s \psi$ are in $\Ve_{a,\e}$. We make use of the fact that $\tilde s$ interchanges the products and that the convolution product is diagonal to represent fusion in terms of a diagonal operator.

\subsection{Fusion in $\Ve_{\e}$} The following theorems allow us to ''diagonalize'' the fusion rules for a modular $\Z_2$ equivariant fusion category. By now we are accustomed to the result in the case of the standard Verlinde algebra. In the general theory it is known that one cannot find such a diagonalization. The grading and action that $\Z_2$ can impose on our category is indeed somewhat limited. Fusion in this algebra is commutative while in the general theory it is not. Regardless- we do not have {\it a priori} that fusion can be diagonalized. In later sections we shall push these results to groups of prime order.

\begin{theorem}\label{t:z2diag} Let $\CC$ be a modular $\mathbb{Z}_{2}$-equivariant fusion category. Fix a generic basis for $\Ve(\CC)$ as described above. Fix $i \in I_{a}$. Let $L_{i}$ be the operator of left multiplication by $\lambda_{i}$ in our basis. Let $M$ be the change of basis operator described in \seref{s:diag} which makes convolution diagonal. Define a linear operator $D_{i} \colon \Ve_{*,\e} \to \Ve_{*,\e}$ by the following equations. For $k \in {}^aI_{\e}$ define

\begin{gather*}
  D_{i} \alpha_{k} = -\frac{(s \lambda_{i}, \basis[a]{k})}{(s \lambda_{0}, \lambda_{k})} \alpha_{k}, \mspace{18mu} D_{i} \beta_{k}  = \frac{(s \lambda_{i}, \basis[a]{k})}{(s \lambda_{0}, \lambda_{k})} \beta_{k}, \\
  \intertext{and for $k \notin {}^aI_{\e}$ and $k \in I_{\e}$ define} 
  D_{i} \lambda_{k} = \frac{(s \lambda_{i}, \basis[a]{k})}{(s \lambda_{0}, \lambda_{k})} \lambda_{k}.
\end{gather*}

Then we have that $\mspace{9mu} M s L_{i} \mspace{3mu} = \mspace{3mu} D_{i} M s$.
\end{theorem}

\begin{proof}
Multiply both sides by $D$ and show that $\mspace{3mu} M \tilde s L_{i} \mspace{3mu} = \mspace{3mu} D_i M \tilde s$. Let $\lambda_{j} \in \Ve_{\e,a}$. Then from the right-hand side of the equation we have $D_{i} M \tilde s \lambda_{j}$

  \[
  \begin{split}
  =& D_{i} M \sum_{k \, \in \, {}^aI_{\e}} \mspace{-9mu} d_k \mspace{9mu}
    \raisebox{-0.5\height}{\begin{overpic}
      {smatrix}
      \put(7,49){$\lambda$}
      \put(3,20){$V_j$}
      \put(-2,65){$V_j$}
      \put(65,5){$V_k$}
      \put(65,90){${}^{a}V_k$}
    \end{overpic}} \\
  =& D_{i} M \sum_{k \, \in \, {}^aI_{\e}} \mspace{-9mu} (\tilde s 
    \lambda_{j}, \basis[a]{k}) \basis[a]{k} \\
  =& D_{i} \sum_{k \, \in \, {}^aI_{\e}} \mspace{-9mu} (\tilde s 
  \lambda_{j}, \basis[a]{k}) (\beta_{k} + \alpha_{k}) \\
  =& \sum_{k \, \in \, {}^aI_{\e}} \mspace{-9mu} \frac{(\tilde s \lambda_{j}
    , \basis[a]{k})(\tilde s \lambda_{i}, \basis[a]{k})}{(\tilde s \lambda_{0},
    \lambda_{k})} (\beta_{k} - \alpha_{k}).
  \end{split}
  \]

Consider the left-hand side. We have that $\tilde s L_{i} \lambda_{j} = \tilde s (\lambda_{i} \ttt \lambda_{j}) = (\tilde s \lambda_{j} * \tilde s \lambda_{i})$

  \[
  \begin{split}
  =& \sum_{k \, \in \, \vphantom{I}^{a}I_{\e}} \mspace{-3mu} (\tilde s \lambda_j)_k \mspace{9mu} * \mspace{9mu} \sum_{k \, \in \, \vphantom{I}^{a}I_{\e}} \mspace{-3mu} (\tilde s \lambda_i)_k \\
  =& \sum_{k \, \in \, \vphantom{I}^{a}I_{\e}} \mspace{-9mu} d_k \mspace{9mu}
      \raisebox{-0.5\height}{\begin{overpic}
      {smatrix}
      \put(7,49){$\lambda$}
      \put(3,20){$V_j$}
      \put(-2,65){$V_j$}
      \put(65,5){$V_k$}
      \put(65,90){${}^{a}V_k$}
    \end{overpic}} \mspace{27mu} * 
    \mspace{9mu} \sum_{k \, \in \, \vphantom{I}^{a}I_{\e}} \mspace{-9mu} d_k \mspace{9mu}
      \raisebox{-0.5\height}{\begin{overpic}
      {smatrix}
      \put(7,49){$\lambda$}
      \put(3,20){$V_i$}
      \put(-2,65){$V_i$}
      \put(65,5){$V_k$}
      \put(65,90){${}^{a}V_k$}
    \end{overpic}} \mspace{27mu} \\ 
  =& \sum_{k \, \in \, \vphantom{I}^{a}I_{\e}} \mspace{-9mu} d_k \mspace{9mu}
      \raisebox{-0.5\height}{\begin{overpic}
      {smatrix}
      \put(7,49){$\lambda$}
      \put(3,20){$V_j$}
      \put(-2,65){$V_j$}
      \put(65,5){$V_k$}
      \put(65,90){${}^{a}V_k$}
    \end{overpic}} \mspace{27mu} * \mspace{9mu} d_k 
    \mspace{9mu}
      \raisebox{-0.5\height}{\begin{overpic}
      {smatrix}
      \put(7,49){$\lambda$}
      \put(3,20){$V_i$}
      \put(-2,65){$V_i$}
      \put(65,5){$V_k$}
      \put(65,90){${}^{a}V_k$}
    \end{overpic}} \mspace{27mu} \\
    =& \sum_{k \, \in \, {}^{a}I_{\e}} \mspace{-9mu} (\tilde s \lambda_{j},
       \basis[a]{k}) \raisebox{-0.5\height}{\begin{overpic}
    {arrow}
    \put(10,5){$V_{k}$}
    \put(10,90){${}^{a}V_{k}$}
  \end{overpic}} * (\tilde s \lambda_{i}, \basis[a]{k}) \raisebox{-0.5\height}{\begin{overpic}
    {arrow}
    \put(10,5){$V_{k}$}
    \put(10,90){${}^{a}V_{k}$}
  \end{overpic}} \\
    =& \sum_{k \, \in \, {}^{a}I_{\e}} \mspace{-9mu} \frac{(\tilde s \lambda_{j}
    , \basis[a]{k})(\tilde s \lambda_{i}, \basis[a]{k})}{(\tilde s \lambda_{0},
       \lambda_{k})} \lambda_{k}.
  \end{split}
  \]

Now apply $M$. This yields

  \[
   \sum_{k \, \in \, {}^{a}I_{\e}} \mspace{-9mu} \frac{(\tilde s \lambda_{j}
    , \basis[a]{k})(\tilde s \lambda_{i}, \basis[a]{k})}{(\tilde s \lambda_{0},
    \lambda_{k})} (\beta_{k} - \alpha_{k})
  \]

and the result.

\end{proof}

As done twice previously we find a formula for the fusion coefficients. The proof runs parallel to \thref{th:verform} and \coref{c:v11mod} and makes use of the change of basis operator. 

\begin{corollary}\label{c:v1fu}

Keep the hypotheses from the above lemma. Then the fusion coefficients are given by the following formula.

  \[
  L_{ij}^k = \sum_{p} \frac{(s \lambda_i, \basis[a]{p})(s \lambda_j, \basis[a]{p})(s \lambda_k^*, \lambda_{p})}{(s \lambda_0, \lambda_p)}. 
  \]
  
\end{corollary}  

\begin{proof}

  Fix $\lambda_{i}$ and $\lambda_{j}$ in $\Ve_{\e,a}$. Then the equation $s L_i = M^{-1} D_{i} M s$ can be written as follows.

 \begin{align*}
   &\mspace{54mu} & s L_i \lambda_j &= M^{-1} D_i M s \lambda_j&\\
   \iff&\mspace{54mu} & s \sum_{r \in I_{\e}} L_{ij}^r \lambda_r &= M^{-1} D_i M \sum_{p \in {}^{a}I_{\e}} (s \lambda_j, \basis[a]{p}) \basis[a]{p} &\\
   \iff&\mspace{54mu} &\sum_{r} L_{ij}^r s \lambda_r &= \sum_{p} M^{-1} D_i M (s \lambda_j, \basis[a]{p}) \basis[a]{p} &\\
   \iff&\mspace{54mu} &\sum_{r} L_{ij}^r s \lambda_r &= \sum_{p} M^{-1} D_i (s \lambda_j, \basis[a]{p}) (\beta_{p} + \alpha_{p}) &\\ 
   \iff&\mspace{54mu} &\sum_{r} L_{ij}^r \sum_{p} (s \lambda_r, \lambda_{p}) \lambda_p &= \sum_{p} M^{-1} \frac{(s \lambda_i, \basis[a]{p})(s \lambda_j, \basis[a]{p})}{(s \lambda_0, \lambda_p)} (\beta_{p} - \alpha_{p}) &\\
   \iff&\mspace{54mu} &\sum_{r} L_{ij}^r \sum_{p} (s \lambda_r, \lambda_{p}) \lambda_p &= \sum_{p} \frac{(s \lambda_i, \basis[a]{p})(s \lambda_j, \basis[a]{p})}{(s \lambda_0, \lambda_p)} \lambda_{p} &\\
   \iff&\mspace{54mu} &\sum_{r} L_{ij}^r (s \lambda_r, \lambda_{p}) &= \frac{(s \lambda_i, \basis[a]{p})(s \lambda_j, \basis[a]{p})}{(s \lambda_0, \lambda_p)}.& 
\end{align*}

We multiply both sides by $(s \lambda_{p}, \lambda_k^*)$ and sum over $p$. Since $s$ is a symmetric unitary operator we have that $(s \lambda_r, \lambda_{p})(s \lambda_{p}, \lambda_k^*) = (s \lambda_r, \lambda_{p})(s \lambda_{p}^*, \lambda_k) = \delta_{rk}$ and the proof.

\end{proof}

These results predict the remainder of the fusion table in $\Ve_{\e}$ in the example from \seref{s:uqsl}. 

\subsection{Example: Fusion in Algebra of Type $D_{2m+2}$}

Consider again the example discussed in \seref{s:uqsl}. Here we want to show that \coref{c:v1fu} correctly predicts the fusion decomposition for $\lambda_{i} \ttt \lambda_{j}$ where $\lambda_{i}, \lambda_{j} \in \Ve_{\e,a}$. The argument is the same as it was in \seref{s:uqsl}. We need to prove \thref{t:twosums} again except here $i$ is odd rather than even. There is nothing to do but change the left-hand side of the equation from \coref{c:v11mod} to \coref{c:v1fu}. As we did in \seref{s:uqsl}, use \thref{t:coefrelat}, which relates the fusion coefficients of the two algebras, and we have shown that \coref{c:v1fu} gives the fusion rules.

\newpage	

\end{document}